\documentclass[12pt]{article}
\usepackage{fullpage,graphicx,amsmath,amsfonts,amsthm}
\usepackage{bm,enumitem,subfig,multirow}
\usepackage{algorithm, algpseudocode}
\usepackage{lineno,hyperref}
\newtheorem{definition}{Definition}[section]
\newtheorem{lemma}{Lemma}[section]
\newtheorem{theorem}{Theorem}[section]
\DeclareMathOperator{\diag}{\text{diag}}

\title{Bounds and Convex Heuristics for Bi-Objective {Optimal Experiment Design} in Water Networks}
\author{Filippo Pecci, Ivan Stoianov}
\date{%
    {\small InfraSense Labs, Department of Civil and Environmental Engineering, Imperial College London, London, SW7 2BB, UK}\\[2ex]%
}

\begin{document}
\maketitle

\begin{abstract}
{Optimal Experiment Design for parameter estimation in water networks} has been traditionally formulated to maximize either hydraulic model accuracy or spatial coverage. Because a unique sensor configuration that optimizes both objectives may not exist, these approaches inevitably result in sub-optimal configurations with respect to one of the objectives.
This paper presents a new bi-objective optimization problem formulation to investigate the trade-offs between these conflicting objectives. We develop a convex heuristic to approximate the Pareto front, and compute guaranteed bounds to discard portions of the criterion space that do not contain non-dominated solutions. Our method relies on a Chebyshev scalarization scheme and convex optimization. We implement the proposed methods for {optimal experiment design} in an operational water network from the UK. For this case study, the convex heuristic computes near-optimal solutions
for the individual objective minimization problems, and tight bounds on the
true Pareto front of the considered bi-objective optimization problem.
\end{abstract}

\section{Introduction}
Hydraulic models play a critical role in supporting water utilities achieving their operational objectives. These include optimal pressure control~\cite{Wright2015}, detection and localization of water leaks~\cite{Blocher2021,Steffelbauer2022}, optimal network design~\cite{Bragalli2015,Ulusoy2021}, improvement of energy efficiency~\cite{Jung2015,Menke2016}.
Monitoring hydraulic pressure and flows with high temporal and spatial resolution is critical to build and maintain accurate hydraulic models~\cite{Zhou2018,Waldron2020}. In addition, monitoring devices are used to provide confidence in the implementation of advanced control solutions.

The problem of {selecting the best sensor locations for parameter estimation in hydraulic models} {belongs to the framework of Optimal Experiment Design (OED)~\cite{Pukelsheim2006,Joshi2009,Kekatos2012,Yu2021,Allen-Zhu2021}}. This problem has traditionally been formulated with the objective of minimizing the uncertainty associated with the hydraulic model, estimated using model sensitivities with respect to pipe roughness parameters~\cite{Bush1998,Kapelan2003,Kapelan2005,Behzadian2009}. However, since the most sensitive locations can be close to each other, such approach may result in unsatisfactory locations for wider monitoring purposes, with sensors clustered around few network nodes - see also the numerical experiments reported in this manuscript. Hence, researchers have considered the problem of {optimal experiment design} in water networks to maximize spatial coverage, using metrics based on network topology~\cite{DeSchaetzen2000,Fontana2015}.

Because a single sensor configuration that is optimal for both sensitivity and topology based metrics may not exist, we propose a novel bi-objective problem formulation for {optimal experiment design} in water networks. As common in {previous literature \cite{Joshi2009,Kekatos2012,Yu2021,Allen-Zhu2021}}, we consider D-optimality as sensitivity-based metric to evaluate the performance of a given sensor configuration. Such convex non-linear objective function aims to minimize the volume of the confidence ellipsoid around the estimated hydraulic parameters. In comparison, the topology metric used by \cite{Fontana2015} aims to generate a spatially even distribution of sensors, by minimizing the sum of the distances between network nodes and their closest sensors. We formulate a bi-objective problem to optimize D-optimality and the topology-based metric. Furthermore, the placement of sensors at network nodes is modeled by binary decision variables, yielding a convex bi-objective mixed integer non-linear optimization problem.

In bi-objective optimization, the best compromises between the two objective functions are characterized by the notion of Pareto optimality.
The Pareto optimal set is made of those feasible solutions such that no other solution exist that reduces one objective function value without increasing the other. The image of the Pareto optimal set through the objective functions is called Pareto front, and it is included in the criterion space, i.e. the space of the objective function values. Due to the presence of binary decision variables, the Pareto front of the considered problem may be disconnected and composed of isolated points. This complicates the generation of the complete Pareto front, and exact solution methods rely on computationally intensive techniques such as branch and bound \cite{DeSantis2020,DeSantis2020a}. 
As a result, the considered bi-objective {optimal experiment design} problems can become computationally impractical to solve exactly when large water networks are considered. Here, we focus on techniques to approximate the Pareto front, discarding portions of criterion space that do not contain any point of the true Pareto front. 

We develop a Chebyshev scalarization scheme, which converts the original bi-objective optimization problem into a series of parameterized, single-objective, mixed integer programs (MIPs). These MIPs are mixed-integer semidefinite programs whose size grows with $n_n^2$, where $n_n$ is the number of nodes in the network. Hence, the computation of certified globally optimal solutions can be impractical when large water networks are considered. We implement a heuristic based on convex optimization to generate feasible solutions with lower bounds on the globally optimal values. 
Furthermore, we investigate the use of feasible solutions and their associate lower bounds to derive an approximation of the Pareto Front with guaranteed bounds. In particular, we prove that the lower bounds of the Chebyshev scalarized MIPs define a subset of the criterion space that does not contain the image of any feasible solution through the objective functions. We evaluate the proposed methods to compute the Pareto front for {optimal experiment design} in an operational water network from the UK. An additional case study with analogous results in presented in the Supplementary Material S4.


\section{Background and problem formulation}
\label{sec:probform}
We consider a water distribution network with $n_p$ links, $n_n$ {demand} nodes, {and $n_0$ known-head nodes (i.e. water inlets)} operating over a time interval $\{1,\ldots,n_t\}$. Let ${h_{0k}} \in \mathbb{R}^{n_0}$ be a known vector of hydraulic heads at water inlets, $d_k \in \mathbb{R}^{n_n}$ be a known vector of nodal demands, and $\eta_k \in \mathbb{R}^{n_v}$ be a known vector of head losses {at the $n_v$ pressure control valves}, for all $k \in \{1,\ldots,n_t\}$. Let $\theta \in \mathbb{R}^{n_r}$ be a vector of roughness parameters, where each index $r \in \{1,\ldots,n_r\}$ corresponds to a group of pipes that share the same roughness parameter $\theta_r$~\cite{Mallick2002}. {We denote with $A_{12}\in\mathbb{R}^{n_p\times n_n}$, $A_{10} \in \mathbb{R}^{n_p \times n_0}$ and $A_{13} \in \mathbb{R}^{n_p \times n_v}$ the link-demand node, link-known head node and link-valve incidence matrices, respectively.  At each time step $k \in \{1,\ldots,n_t\}$}, network hydraulic heads $h_k \in \mathbb{R}^{n_n}$ and flows $q_k \in \mathbb{R}^{n_p}$ are computed by solving the following non-linear equations:
\begin{equation}
\label{eq:hyd_eqs}
    \begin{split}
        &A_{11}(q_k,\theta)q_k + A_{12} h_k + A_{10}{h_{0k}} + A_{13}\eta_k = 0\\
        &A_{12}^Tq_k= d_k.
    \end{split}
\end{equation}
The vector $A_{11}(q_k,\theta)q_k$ corresponds to the head losses within network links due to friction, and the diagonal matrix $A_{11}(q_k,\theta)$ is defined by:
\begin{equation}
    A_{11}(q_{lk},\theta)_{ll} = \rho_{l}(q_{lk},\theta_{r_l})|q_{lk}|^{n_l-1}, \quad l \in \{1,\ldots,n_p\},
\end{equation}
where $r_l \in \{1,\ldots,n_r\}$ is the index of the group corresponding to link $l \in \{1,\ldots,n_p\}$, and $\rho_l$ and $n_l$ are resistance coefficient function and exponent associated with link $l$, respectively, whose formulations depend on the head loss formula used \cite[Section 2.2]{Larock1999}. In the case of Hazen-Williams friction model, we have $n_l=1.852$ and
\begin{equation}
\rho_l(\theta_{r_l}) = \frac{10.67 L_l}{\theta_{r_l}^{1.852} D_l^{4.871}},
\end{equation}
{where $L_l$ and $D_l$ are length and diameter of link $l$, respectively.} When the Darcy-Weisbach formula is used, we have $n_l=2$ and $\rho_l(q_{lk},\theta_{r_l})$ is a continuously differentiable function whose piece-wise formulation can be found in \cite{Waldron2020}.

In the following, we denote by $\phi_k:\mathbb{R}^{n_r} \rightarrow \mathbb{R}^{n_p+n_n}$ the non-linear function such that $\phi_k(\theta) = [q_k^T\;h_k^T]^T$ solves \eqref{eq:hyd_eqs}. 
Let $S \subset \{1,\ldots,n_n\}$ and $S_0 \subset \{1,\ldots,n_p\}$ be index sets corresponding to hydraulic heads and flow measuring devices, respectively. Let measurement $y_{ik}$ correspond to either a measured flow or a measured hydraulic head at time step $k$, for each $i \in S \cup S_0$.

As done in previous literature, we consider pipe roughness coefficients $\theta \in \mathbb{R}^{n_r}$ as parameters to be estimated. 
It is possible to estimate the vector of parameters $\theta \in \mathbb{R}^{n_r}$ by solving a non-linear inverse problem:
\begin{equation}
\label{eq:param_est_prob}
\begin{alignedat}{3}
&\text{minimize}&\; \; & \frac{1}{n_t|S \cup S_0|}\sum_{k=1}^{n_t}\sum_{i \in S \cup S_0}\mathcal{L}(\phi_{ik}(\theta),y_{ik}) + g(\theta)\\
&\text{subject to}& & \theta \geq 0
\end{alignedat}
\end{equation}
where $\mathcal{L}(\cdot,y_{ik})$ and $g$ are suitable convex loss and regularization functions, respectively. For examples on parameter estimation in water networks see \cite{Lansey1991} and \cite{Waldron2020}.

The quality of the parameter estimates obtained solving{~\eqref{eq:param_est_prob}} heavily depends on the choice of the measuring locations $S$ and $S_0$. 
The problem of selecting the best sensor locations for parameter estimation has been formulated within the framework of {Optimal Experiment Design (OED)~\cite{Pukelsheim2006}}. We assume that the vector of parameters $\theta \in \mathbb{R}^{n_r}$ and the measurements are related through the following Gaussian additive non-linear model:
\begin{equation}
    \label{eq:inv_prob}
    y_{ik} = \phi_{ik}(\theta) + v_{ik}, \quad i\in S \cup S_0
\end{equation}
with $v_{ik} \sim N(0,\sigma^2)$ independent identically distributed random variables. The goal of {OED} is to select sensor locations so that the estimation error is minimized.
In the case of \eqref{eq:inv_prob}, the error covariance matrix is computed using the sensitivity vectors:
\begin{equation}
   \Sigma = \sigma^2\bigg(\sum_{j \in S_0}A_jA_j^T + \sum_{j \in S}A_jA_j^T  \bigg)^{-1},
\end{equation}
with $A_j = [\nabla \phi_{j,1}(\bar{\theta}) \ldots \nabla \phi_{j,n_t}(\bar{\theta})] \in \mathbb{R}^{n_r \times n_t}$, given an initial estimate of the parameter vector denoted by $\bar{\theta} \in \mathbb{R}^{n_r}$. {Recall that $A_j$ will either refer to a flow measurement if $j \in S_0$, or to a pressure measurement if $j \in S$.} 
In this study, we exploit the saddle-point structure of Equations~\eqref{eq:hyd_eqs} and derive a computationally efficient algorithm to compute these sensitivity vectors - see supplementary material S1 for a detailed description. 
We assume that $n_r < |S \cup S_0|$ and vectors $\nabla \phi_{jk}(\bar{\theta})$ span $\mathbb{R}^{n_r}$. If this is not the case, a prior should be imposed on the parameters $\theta \in \mathbb{R}^{n_r}$. If $\theta_l \sim \mathcal{N}(0,\lambda)$, we obtain a covariance matrix $\tilde{\Sigma} = (\frac{1}{\sigma^2}\sum_{j \in S \cup S_0}A_jA_j^T + \frac{1}{\lambda} I_r)^{-1}$, where $I_r \in \mathbb{R}^{n_r \times n_r}$ is the identity matrix {- see Equation~(16) in \cite{Joshi2009}}. Formulations and methods discussed in the remainder of this manuscript are then applied to $\tilde{\Sigma}$. In addition, we assume that the locations of flow measuring devices $S_0 \subset \{1,\ldots,n_p\}$ are fixed, and we want to optimize the choice of $S \subset \{1,\ldots,n_n\}$. Let binary variables $z \in \{0,1\}^{n_n}$ be used to model the placement of sensors, with $z_j=1$ if a sensor is placed on node $j$ or $z_j=0$ otherwise. Moreover, we set $X_0 = \sum_{j \in S_0}A_jA_j^T$. 

Several metrics have been proposed to measure the estimation error associated with the covariance matrix $\Sigma$. A popular choice is D-optimality~\cite{Bush1998,Joshi2009,Kekatos2012,Yu2021, Allen-Zhu2021}, which aims to minimize the volume of the confidence ellipsoid around the estimated parameters:
\begin{equation}
    f_D(z) = -\log\,\det\bigg(X_0+\sum_{j=1}^{n_n}z_jA_jA_j^T\bigg).
\end{equation}

Metrics used in {OED literature \cite{Kapelan2003,Kapelan2005,Joshi2009,Behzadian2009,Kekatos2012,Yu2021,Allen-Zhu2021}}  rely on sensitivity information to evaluate the estimation error. However, they do not explicitly consider network topological properties. Since the most sensitive locations can be very close to each other, sensitivity-based metrics may result in unsatisfactory sensor configurations. For example, when the number of installed sensors is greater than $2n_r$, it is possible to place sensors at both inlet and outlet of a link for each group $r \in \{1,\ldots,n_r\}$. In this way, the head loss across the link would be measured, reducing the uncertainty associated with the estimated roughness parameter for the corresponding group. For these reasons, researchers have also considered topology-based methods, which do not require an initial estimate of the parameters $\bar{\theta}$, and take into account network connectivity to generate an even spread of sensors~\cite{DeSchaetzen2000}. In particular, \cite{Fontana2015} proposed a facility location problem to minimize: 
\begin{equation}
    f_T(z) = \sum_{i=1}^{n_n} \min_{j \in N_i \cap I(z)} C_{ij},
\end{equation}
where $N_i = \{1,\ldots,n_n\}\setminus \{i\}$ and $I(z) = \{ j \in \{1,\ldots,n_n\} \; | \; z_j=1\}$. Coefficients $C_{ij}$ are defined as the length of the shortest-path between node $i$ and node $j$, where each link $l$ has been assigned a weight $L_l$, with $l \in \{1,\ldots,n_p\}$. Function $f_T$ can be interpreted as the sum of the distances between network nodes and their closest sensors. The following Lemma shows how to reformulate $f_T$ using linear programming.
\begin{lemma}
\label{lem:pmedian_equivalently}
Let $\hat{z} \in \{0,1\}^{n_n}$, then $f_T(\hat{z})$ is the optimal value of the following linear program:
\begin{equation}
    \label{eq:pmedian_z}
    \begin{alignedat}{3}
    &\text{minimize}&\; \; & \sum_{i=1}^{n_n}\sum_{j=1}^{n_n} C_{ij}u_{ij}\\
    &\text{subject to}& & u_{ij} \leq \hat{z}_j,\quad i \in \{1,\ldots,n_n\},\; j \in \{1,\ldots,n_n\}\\
    &&&\sum_{j=1}^{n_n}u_{ij} = 1, \quad  i \in \{1,\ldots,n_n\}\\
    &&& u_{ii} = 0, \quad  i \in \{1,\ldots,n_n\}\\
    &&& u_{ij} \geq 0 \quad i \in \{1,\ldots,n_n\},\; j \in \{1,\ldots,n_n\}
    \end{alignedat}
\end{equation}
\end{lemma}
\begin{proof}
    See section S2 of the Supplementary Material.
\end{proof}
The above lemma implies that the problem of minimizing $f_T$ can be formulated as a mixed integer linear program (MILP) and solved using integer programming techniques~\cite{Fontana2015}. Because the size of these MILPs grows rapidly with the number of nodes, previous literature has also considered both specialized heuristics~\cite{Reese2006,Ushakov2021}, and model reduction approaches for specific network topologies~\cite{Church2020}. 

As shown by the numerical experiments reported in Section~\ref{sec:numexp}, a single sensor configuration that simultaneously minimizes $f_D$ and $f_T$ may not exist. Therefore, in the next Section, we introduce a new bi-objective problem formulation, and a scalarization method to generate an approximation of the associated Pareto front. 
\section{Bi-objective sensor placement in water networks}
\label{sec:multiobj}
To investigate the trade-offs between sensitivity and topology based metrics, we propose a new bi-objective problem formulation:
\begin{equation}
    \label{eq:prob_multiobj}
    \begin{alignedat}{3}
    &\text{minimize}&\; \; &(f_D(z),f_T(z))\\
    &\text{subject to}& & \sum_{j=1}^{n_n}z_j=m\\
    &&& Gz \leq b \\
    &&& z \in \{0,1\}^{n_n},
    \end{alignedat}
\end{equation}
where $m\geq 2$ corresponds to the number of sensors to be installed. We have included constraints to prevent two sensors from being placed on nodes that are connected by a link, avoiding impractical sensor configurations. To achieve this, we define matrix $G \in \mathbb{R}^{n_p \times n_n}$ as:
\begin{equation}
    G(l,j) = \begin{cases} 1& \text{if link $l$ is connected to node $j$}\\
0 & \text{otherwise},
    \end{cases}
    \end{equation}
and vector $b = [1,\ldots,1]^T$. Solving Problem~\eqref{eq:prob_multiobj} aims to compute the set of Pareto optimal solutions, also known as the efficient set. In fact, the most efficient compromises between $f_D$ and $f_T$ are characterized by the notion of Pareto optimality.

\begin{definition}
    A sensor configuration $z^*$ is said to be Pareto optimal for Problem~\eqref{eq:prob_multiobj} if there is no 
    $z'\in \{0,1\}^{n_n}$ feasible for Problem~\eqref{eq:prob_multiobj}, such that $f_U(z') \leq f_U(z^*)$ for all $U \in \{D,T\}$, and $f_Y(z')<f_Y(z^*)$ for some $Y \in \{D,T\}$.
\end{definition}
Let $\Theta \subset \{0,1\}^{n_n}$ be the set of Pareto optimal solutions for Problem~\eqref{eq:prob_multiobj}. Due to the presence of binary decision variables, $\Theta$ is disconnected and composed of isolated points. Because of such complicated structure, researchers have focused on approximating the Pareto front, which is the image of $\Theta$ through the objective functions:
\begin{equation}
    F(\Theta) = \bigg\{\big(f_D(z^*),f_T(z^*)\big)\; | \; z^* \in \Theta \bigg\}.
\end{equation} 
To the best of the authors' knowledge, previous studies have solved {optimal experiment design} problems in water distribution networks using genetic algorithms~\cite{Kapelan2003,Kapelan2005,Behzadian2009}. While being practical approaches for highly non-linear, discrete, or black-box optimization problems, these methods do not provide theoretical guarantee of global non-dominance, and they require a significant computational effort to generate good quality solutions. In comparison, we investigate the application of mathematical optimization methods to approximate the Pareto front with guaranteed bounds. To complete our numerical investigation, Section~\ref{sec:nsgaii} reports a comparison with { three off-the-shelf Multi-objective Optimization Evolutionary Algorithms (MOEAs).} 

Classical algorithms to approximate the Pareto front (known as \emph{criterion space algorithms}) rely on scalarization methods, which convert a bi-objective problem into a series of parameterized single-objective optimization problems. In this study, we investigate a Chebyshev scalarization method \cite[Section 3.4]{Miettinen1999}. As initial step we consider the individual minimization problems:
\begin{equation}
    \label{eq:prob_indiv}
    \begin{alignedat}{3}
    &\text{minimize}&\; \; &f_Y(z)\\
    &\text{subject to}& & \sum_{j=1}^{n_n}z_j=m\\
    &&& Gz \leq b \\
    &&& z \in \{0,1\}^{n_n},
    \end{alignedat}
\end{equation}
for all $Y \in \{D,T\}$. Here, we assume that we have implemented a method to compute a feasible solution $\hat{z}_Y$ and a lower bound to the optimal value of Problem~\eqref{eq:prob_indiv}, denoted by $f^*_Y$, for all $Y \in \{D,T\}$. Section~\ref{sec:singleopt} presents a suitable heuristic based on convex optimization to achieve this. Given $\beta \in (0,1)$, consider the following optimization problem:
\begin{equation}
    \label{eq:prob_cheb}
    \begin{alignedat}{3}
    &\text{minimize}&\; \; &f_{\beta}(z)\\
    &\text{subject to}& & \sum_{j=1}^{n_n}z_j=m\\
    &&& Gz \leq b\\
    &&& z \in \{0,1\}^{n_n},
    \end{alignedat}
\end{equation}
where 
\begin{equation}
    f_{\beta}(z) = \max \bigg(w^{\beta}_D(f_D(z)-f_D^*),w^{\beta}_T(f_T(z)-f^*_T)\bigg),
\end{equation} 
with
\begin{equation}
    w^{\beta}_D = \frac{\beta}{f_D(\hat{z}_T)-f_D^*}, \quad w^{\beta}_T = \frac{1-\beta}{f_T(\hat{z}_D)-f_T^*}. 
\end{equation}
While other scalarization approaches (e.g. weighted sum) are limited to the generation of the convex portions of the Pareto front, any Pareto optimal solution can be found solving Problem~\eqref{eq:prob_multiobj} for an opportunely selected weight. 
{Before stating these results, we introduce the definition of weak Pareto optimality.
\begin{definition}
    A sensor configuration $z^*$ is said to be weakly Pareto optimal for Problem~\eqref{eq:prob_multiobj} if there is no 
    $z'\in \{0,1\}^{n_n}$ feasible for Problem~\eqref{eq:prob_multiobj}, such that $f_U(z') < f_U(z^*)$ for all $U \in \{D,T\}$.
\end{definition}
It follows from the definitions that a Pareto optimal solution is also a weakly Pareto optimal solution.
}
\begin{theorem}
    \label{thm:weakly}
    Let $\hat{z}_{\beta}$ be a globally optimal solution of Problem~\eqref{eq:prob_cheb}. Then, $\hat{z}_{\beta}$ is weakly Pareto optimal.
\end{theorem}
\begin{proof}
    See section S2 of the Supplementary Material.
\end{proof}
\begin{theorem}
    \label{thm:exists}
    Given $z^* \in \Theta$, there exist $\beta^* \in (0,1)$ such that $z^*$ solves Problem~\eqref{eq:prob_cheb} for $\beta = \beta^*$.
\end{theorem}
\begin{proof}
    See section S2 of the Supplementary Material.
\end{proof}
Problem~\eqref{eq:prob_cheb} results in  mixed-integer semidefinite program (SDP) with exponential cone constraints, whose size grows with $n_n^2$ - see \ref{sec:reform}. As a result, obtaining certified globally optimal solutions for Problem~\eqref{eq:prob_cheb} can be impractical when large water networks are considered. For this reason, we only assume that we are able to compute a feasible solution $\hat{z}_{\beta}$ and a lower bound to the optimal value of Problem~\eqref{eq:prob_cheb}, denoted by $f^*_{\beta}$ - Section~\ref{sec:singleopt} presents a suitable convex heuristic. {Let $N>1$ be the number of feasible solutions $\hat{z}_{\beta_1},\ldots,\hat{z}_{\beta_N}$ computed solving Problem~~\eqref{eq:prob_cheb} for different parameters $\beta_1,\ldots,\beta_N$, respectively.} Furthermore, we include the individual minimizers by defining $\hat{z}_{\beta_0} = \hat{z}_D$ and $\hat{z}_{\beta_{N+1}} = \hat{z}_T$. {The number $N$ is set \emph{a priori} based on how many points on the approximated Pareto front we aim to compute}. An approximation of the Pareto front is given by:
\begin{equation}
    \mathcal{P} = \bigg\{(f_D(\hat{z}_{\beta_k}),f_T(\hat{z}_{\beta_k})) \; | \; k=0,\ldots,N+1\bigg\}.
\end{equation}
\begin{theorem}
    \label{thm:upbound}
    Define the following set:
    \begin{equation}
        \label{eq:setdef1}
        \mathcal{W} = \bigg\{(p_1+u,p_2+v) \; | \; (p_1,p_2) \in \mathcal{P},u\geq 0, v \geq 0, u+v>0 \bigg\}.  
    \end{equation}
    Then, we have:
    \begin{equation}
        F(\Theta) \cap \mathcal{W} =  \emptyset 
    \end{equation}
\end{theorem}
\begin{proof}
    Let $z^* \in \Theta$ be such that $(f_D(z^*),f_T(z^*)) \in \mathcal{W}$. Without loss of generality, we assume that $f_D(z^*)=f_D(\hat{z}_{\beta_k})+u$ and $f_T(z^*)=f_T(\hat{z}_{\beta_k})+v$, for some $k \in \{0,\ldots,N+1\}$ and $u \geq 0$, $v >0$. Hence, $f_D(z^*)\geq f_D(\hat{z}_{\beta_k})$ and $f_T(z^*)>f_T(\hat{z}_{\beta_k})$. This contradicts the assumption that $z^* \in \Theta$. Therefore, we must have that $(f_D(z^*),f_T(z^*)) \notin \mathcal{W}$.
\end{proof}
It is also possible to bound the Pareto front from below as follows. Assume to have obtained lower bounds to the optimal value of Problem~\eqref{eq:prob_cheb}, for different values of $\beta$, denoted by $f^*_{\beta_1},\ldots,f^*_{\beta_N}$. Define the following sets:
\begin{equation}
    \label{eq:setdef2}
    \begin{split}
    &\mathcal{L} = \bigg\{\big(\frac{f^*_{\beta_{k}}}{w_D^{\beta_k}}+f^*_D,\, \frac{f^*_{\beta_{k}}}{w_T^{\beta_k}}+f^*_T\big) \; | \; k=1,\ldots,N \bigg\}\\
    &\mathcal{R} = \bigg\{(p_1,p_2) \in \mathbb{R}^2 \; | \; p_1<f_D^* \text{ or } p_2<f_T^* \bigg\}\\
    &\mathcal{Q} = \bigg\{(p_1-u,p_2-v) \; | \; (p_1,p_2) \in \mathcal{L}, u>0, v>0\bigg\}.
    \end{split}
\end{equation}
\begin{theorem}
    \label{thm:bound}
The set $\mathcal{R} \cup \mathcal{Q}$ does not contain the image of any feasible point for Problem~\eqref{eq:prob_multiobj} through the objective functions $(f_D,f_T)$. In particular:
\begin{equation}
    F(\Theta) \cap (\mathcal{R} \cup \mathcal{Q}) = \emptyset
\end{equation}
\end{theorem}
\begin{proof}
    Let $\tilde{z}$ be feasible for Problem~\eqref{eq:prob_multiobj} and assume that $(f_D(\tilde{z}),f_T(\tilde{z}) \in \mathcal{R}$. This implies that either $f_D(\tilde{z})<f_D^*$ or $f_T(\tilde{z})<f_T^*$. But this contradicts the assumption that $f_D^*$ and $f_T^*$ are lower bounds to the optimal values of the individual minimization problems for $f_D$ and $f_T$, respectively. Hence, it must be that $\tilde{z} \notin \mathcal{R}$. Now assume that $(f_D(\tilde{z}),f_T(\tilde{z})) \in \mathcal{Q}$. This implies that there exist $k \in \{1,\ldots,N\}$, $u>0$, and $v>0$ such that
    \begin{equation}
        f_D(\tilde{z}) = \frac{f^*_{\beta_{k}}}{w_D^{\beta_k}}+f^*_D - u, \quad f_T(\tilde{z}) = \frac{f^*_{\beta_{k}}}{w_T^{\beta_k}}+f^*_T -v
    \end{equation}
    Hence, we have:
    \begin{equation}
        f_D(\tilde{z}) < \frac{f^*_{\beta_{k}}}{w_D^{\beta_k}}+f^*_D , \quad f_T(\tilde{z}) < \frac{f^*_{\beta_{k}}}{w_T^{\beta_k}}+f^*_T,
    \end{equation}
    which is equivalent to:
    \begin{equation}
        w_D^{\beta_k}(f_D(\tilde{z})-f^*_D) < f^*_{\beta_{k}} , \quad w_T^{\beta_k}(f_T(\tilde{z})-f^*_T) < f^*_{\beta_{k}}.
    \end{equation}
    This yields $f_{\beta_k}(\tilde{z}) < f^*_{\beta_k}$, which contradicts the assumption that $f^*_{\beta_k}$ is a lower bound to the optimal value of Problem~\eqref{eq:prob_cheb} with $\beta = \beta_k$. Therefore, we conclude that $(f_D(\tilde{z}),f_T(\tilde{z})) \notin \mathcal{Q}$.
\end{proof}

\begin{figure}[h!]
    \subfloat[Case of $N=4$]{\includegraphics[width=0.48\textwidth]{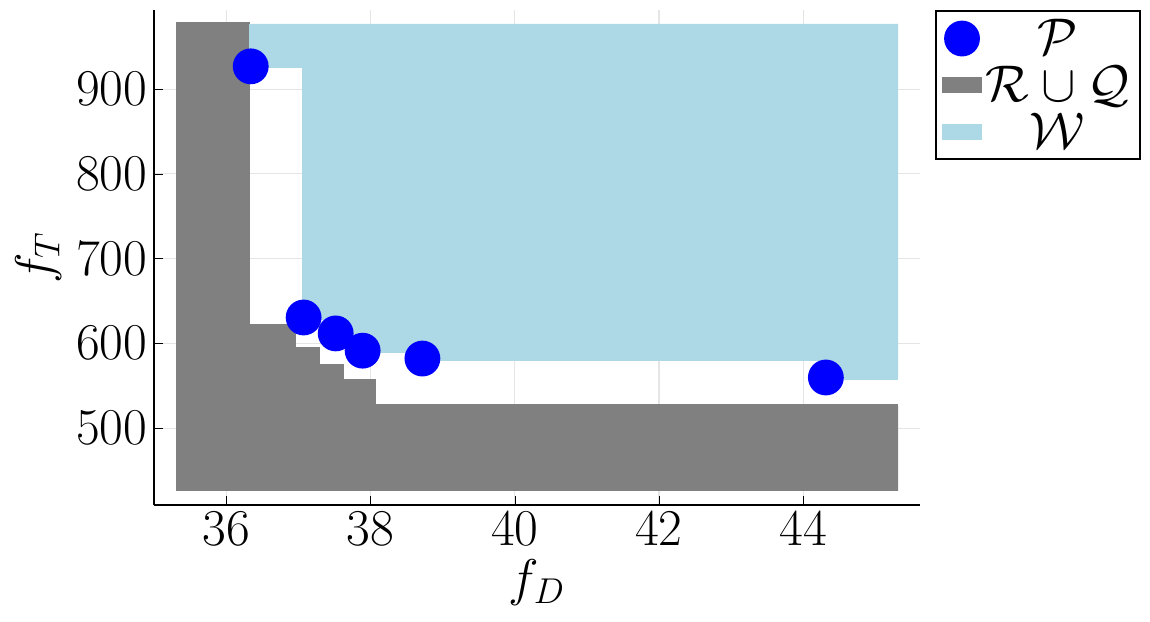}} \, 
    \subfloat[Case of $N=200$]{\includegraphics[width=0.48\textwidth]{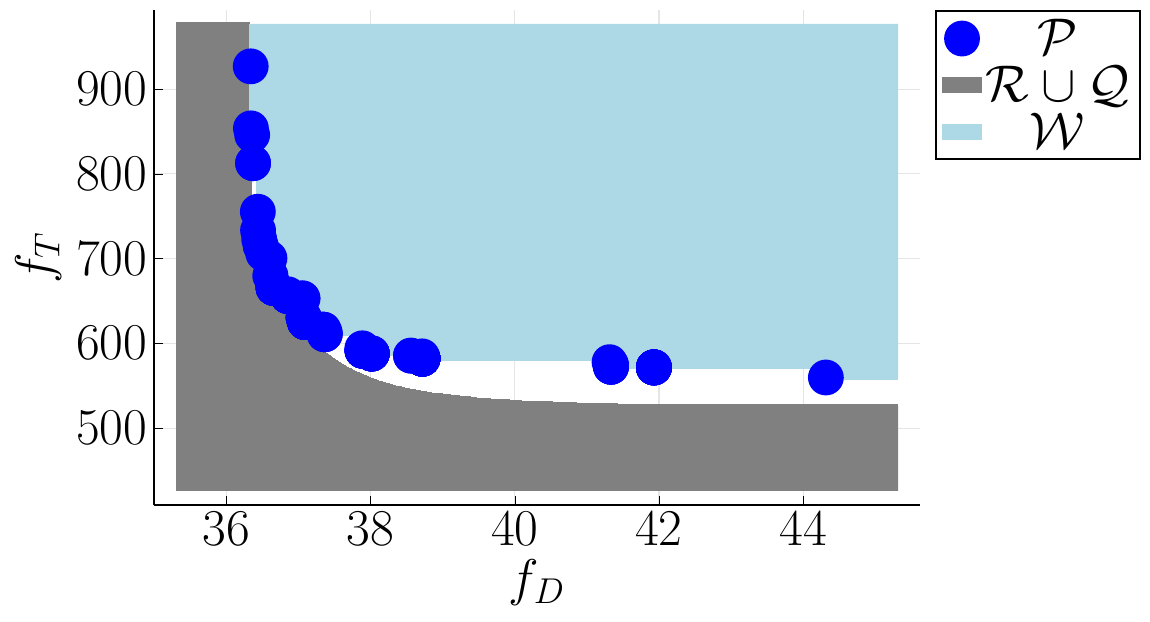}}
    \caption{Examples of Pareto front approximations computed by Algorithm~\ref{alg:cheb}.}
    \label{fig:examples}
\end{figure}

Figure~\ref{fig:examples} shows examples of sets $\mathcal{P}$, $\mathcal{W}$, $\mathcal{R}$, and $\mathcal{Q}$, with different values of $N$. We observe that, as we increase the number of weights, the proposed method can provide a tight outer-approximation of the true Pareto front, which is included in the white area between sets $\mathcal{W}$ and $\mathcal{R}\cup \mathcal{Q}$. 

The proposed Chebyshev scalarization method is summarized in Algorithm~\ref{alg:cheb}. The method returns a Pareto front approximation $\mathcal{P}$ and sets $\mathcal{W}$, $\mathcal{R}$, and $\mathcal{Q}$, which bound the region containing the true Pareto front. 
\begin{algorithm}[h!]
\begin{algorithmic}
    \State{\textbf{Input:} Matrices $A_1$,\ldots,$A_{n_n}$, $C$, and number of scalarization steps $N$.}
    \State{\textbf{Output:} Sets $\mathcal{P}$, $\mathcal{W}$, $\mathcal{R}$, $\mathcal{Q}$.}
    \State{Apply a suitable solver to Problem~\eqref{eq:prob_indiv} to minimize $f_D$ and obtain a feasible solution $\hat{z}_D$ and lower bound $f^*_D$.}
    \State{Apply a suitable solver to Problem~\eqref{eq:prob_indiv} to minimize $f_T$ and obtain a feasible solution $\hat{z}_T$ and lower bound $f^*_T$.}
    \State{Define $\mathcal{P}\gets \{f_D(\hat{z}_D),f_T(\hat{z}_D)\}\cup \{f_D(\hat{z}_T),f_T(\hat{z}_T)\}$ and $\mathcal{L}\gets\emptyset$.}
    \State{Generate a sequence of weights $\beta_1,\ldots,\beta_N$}
    \For{$k=1,\ldots,N$}
    \State{Apply a suitable solver to Problem~\eqref{eq:prob_cheb} with $\beta=\beta_k$ and obtain a feasible solution $\hat{z}_{\beta_k}$ and lower bound $f^*_{\beta_k}$.}
    \State{$\mathcal{P}\gets \mathcal{P} \cup \big\{(f_D(\hat{z}_{\beta_k}),f_T(\hat{z}_{\beta_k}))\big\}$.}
    \State{$\mathcal{L}\gets \mathcal{L} \cup \bigg\{\big(\frac{f^*_{\beta_{k}}}{w_D^{\beta_k}}+f^*_D,\, \frac{f^*_{\beta_{k}}}{w_T^{\beta_k}}+f^*_T\big)\bigg\}$.}
    \EndFor
    \State{Define sets $\mathcal{W}$, $\mathcal{R}$, $\mathcal{Q}$ as in \eqref{eq:setdef1} and \eqref{eq:setdef2}.}
\end{algorithmic}
\caption{Chebyshev scalarization with global bounds.}
\label{alg:cheb}
\end{algorithm}

\section{Solution of single-objective sensor placement problems}
\label{sec:singleopt}
Assuming that either $f(z) = f_D(z)$, $f(z)=f_T(z)$, or $f(z) = f_{\beta}(z)$, Problems~\eqref{eq:prob_indiv} and \eqref{eq:prob_cheb} correspond to the following mixed integer optimization problem:
\begin{equation}
    \label{eq:prob_minlp}
    \begin{alignedat}{3}
    &\text{minimize}&\; \; &f(z)\\
    &\text{subject to}& & \sum_{j=1}^{n_n}z_j=m\\
    &&& Gz \leq b\\
    &&& z \in \{0,1\}^{n_n}.
    \end{alignedat}
\end{equation}
When $f(z)=f_T(z)$, Lemma~\ref{lem:pmedian_equivalently} implies that Problem~\eqref{eq:prob_minlp} is re-written as a mixed integer linear program - see also \ref{sec:reform}. In comparison, when $f(z)=f_D(z)$ and $f(z)=f_{\beta}(z)$, Problem~\eqref{eq:prob_minlp} results in a mixed-integer semidefinite program with exponential cone constraints. In particular, when $f(z)=f_{\beta}(z)$, the size of these mixed-integer problems grows with $n_n^2$. Hence, the {utilization} of exact methods based on branch and bound techniques can become impractical. In this section, we propose a heuristic based on convex optimization to solve Problem~\eqref{eq:prob_minlp}. We develop a convex heuristic based on two steps.
\begin{enumerate}
    \item Relax the binary constraints and solve the following continuous optimization problem to obtain a vector of fractional values $z^* \in [0,1]$ and a lower bound $f^*=f(z^*)$ to the optimal value:
    \begin{equation}
        \label{eq:prob_convex}
        \begin{alignedat}{3}
        &\text{minimize}&\; \; &f(z)\\
        &\text{subject to}& & \sum_{j=1}^{n_n}z_j=m\\
        &&& Gz \leq b\\
        &&& z \in [0,1]^{n_n}.
        \end{alignedat}
    \end{equation}
    \item Implement a heuristic algorithm to generate a good quality feasible solution $\hat{z} \in \{0,1\}^{n_n}$.
\end{enumerate}
Observe that the proposed convex heuristic results in a bound on the level of sub-optimality of $\hat{z}$, given by:
\begin{equation}
    \text{Gap} = \frac{f(\hat{z})-f^*}{f^*}
\end{equation}
When $f(z)=f_T(z)$, Lemma~\eqref{lem:pmedian_equivalently} implies that Problem~\eqref{eq:prob_convex} is a equivalently re-written as a linear program. Hence, it can be efficiently solved by state-of-the-art linear solvers e.g. GUROBI \cite{GurobiOptimization2021}. In the cases of $f_D$ and $f_{\beta}$, we reformulate the convex non-linear Problem~\eqref{eq:prob_convex} as a SDP with exponential cone constraints, and solve it using the convex optimization solver MOSEK~\cite{mosek}. Since the size of the considered SDPs grows with $n_n^2$, it might be necessary to implement tailored decomposition approaches \cite{Garstka2021} to derive scalable solution methods for problems formulated on large water networks - e.g. with $n_n>10^4$. These SDP reformulations are needed because the direct implementation of state-of-the-art non-linear programming solver IPOPT~\cite{Wachter2006} resulted in a challenging computational performance, as discussed in Section~\ref{app:ipopt}. Details on these reformulations are given in \ref{sec:reform}. In the remainder of this Section, we present the heuristic algorithm. 

Let $z^* \in [0,1]$ be a solution of Problem~\eqref{eq:prob_convex}, and $f^*=f(z^*)$ be the corresponding lower bound to the optimal value of Problem~\eqref{eq:prob_minlp}. We investigate a heuristic method to compute a feasible solution for Problem~\eqref{eq:prob_minlp}, when either $f(z)=f_D(z)$, $f(z)=f_T(z)$, or $f(z)=f_{\beta}(z)$. The proposed method builds on a Fedorov exchange algorithm by \cite{Joshi2009}, adding the handling of constraints on allowed sensor locations and different objective functions in addition to $f_D$. We refer to the implemented algorithm as Round-and-Swap, see Supplementary Material S3 for a pseudo-code description. 

As first step, the fractional values in $z^*$ are rounded to obtain a vector of binary values $\hat{z} \in \{0,1\}^{n_n}$ that satisfies the problem's constraints. Since we expect the largest elements in $z^*$ to correspond to good quality sensor configurations, we round the largest $m$ elements in $z^*$ to $1$, and to $0$ the others. We also check that each selected location does not result in violation of the linear constraints, in which case we select the next location according to the descending order of $z^*$-values. 

Next, we investigate alternative sensor configurations by considering those that are obtained swapping one of the $m$ locations identified by $\hat{z}$ for one of the $n_n-m$ locations that were not selected. When we find a feasible sensor configuration that reduces the objective function value, we update $\hat{z} \in \{0,1\}^{n_n}$. If we find that no swap reduces the objective function {value of} the current solution, the algorithm terminates. In comparison to the algorithm in \cite{Joshi2009}, at each swapping iteration, we do not only evaluate the objective function value of the new sensor configuration, but we also check that it satisfies the linear constraints. In this way, the algorithm is guaranteed to generate feasible solutions. 

While the swapping process is guaranteed to terminate after a finite number of iterations, it can still require a large number of function evaluations. Hence, we define a maximum number of iterations $\texttt{MaxIter}$, which is set to $10^5$ in our {numerical experiments}. The largest entries in $z^*$ are expected to include the globally optimal sensor locations, and this can help to reduce the number of iterations needed to converge. In particular, we consider for swapping only locations corresponding to entries in $z^*$ with values between $0.1$ and $0.9$. We identify locations to be removed from the set of selected sensor locations in ascending order based on their $z^*$-values. Analogously, the sensor locations to be added to the set of selected locations are chosen in descending order based on their values in $z^*$.

\section{Numerical experiments}
\label{sec:numexp}
All experiments reported in this Section are carried out in Julia 1.7.2~\cite{Julia-2017}, where optimization solvers are called through JuMP 1.0.0~\cite{DunningHuchetteLubin2017}. All computations are performed on a laptop with Intel(R) Core(TM) i9-9980HK CPU @ 2.40GHz and 32 GB of RAM. Linear programs~\eqref{eq:pmedian_convex} are solved by GUROBI~\cite[version 9.1]{GurobiOptimization2021}. Semidefinite programs \eqref{eq:dopt_convex} and~\eqref{eq:cheb_convex3} are solved using MOSEK~\cite[version 9.3]{mosek}.

\begin{figure}
    \centering
    \subfloat[Network layout. \label{fig:bwflnet}]{\includegraphics[width=0.48\textwidth]{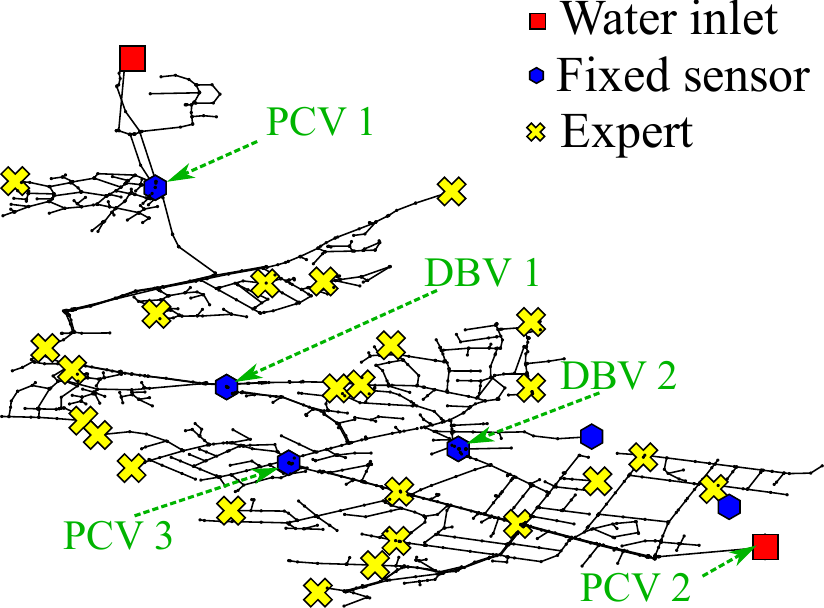}}
    \,
    \subfloat[Pipe groups. \label{fig:bwflnet_groups}]{\includegraphics[width=0.48\textwidth]{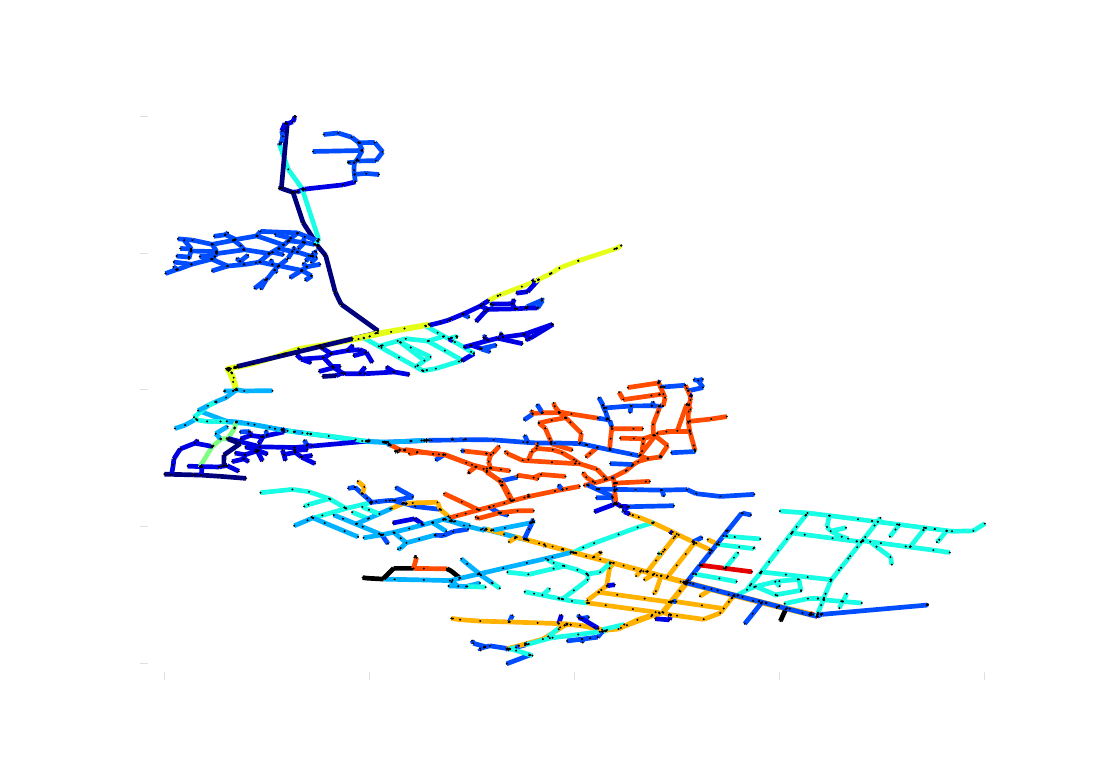}}
    \label{fig:bwflnet_all}
    \caption{\texttt{BWFLNet}: (a) shows locations of fixed sensors, expert choice of existing pressure sensors, pressure control valves (PCVs), dynamic boundary valves (DBVs); (b) pipe groups where black color is used to highlight valves and leaf pipes that are not considered by the {optimal experiment design} problem. }
\end{figure}

We implement the developed methods for optimizing the locations of pressure sensors to be installed in \texttt{BWFLnet}, the hydraulic model of an operational water network from the UK - see Figure~\ref{fig:bwflnet}. The network has $n_0=2$ inlets, $n_n=2221$ nodes, and $n_p=2281$ links, including $3$ pressure control valves (PCVs) which are optimally controlled to minimize Average Zone Pressure~\cite{Wright2015}. The EPANET model for \texttt{BWFLnet} is available at the Mendeley data repository in \cite{Pecci2021}.
In addition, \texttt{BWFLnet} is operated with a dynamically adaptive topology, where two dynamic boundary valves (DBVs) are open during diurnal operation and closed at night for leakage management purposes. The network also includes a number of throttle control valves (TCVs), that are assumed to be fully open. We consider all valves' coefficients as known and we do not perform parameter estimation on these links. The hydraulic model has $13$ different pipe groups, which were determined by the network operator. Two of these groups, namely groups $8$ and $9$, consists of leaf links whose roughness parameters can be estimated only by installing sensors at both ends of the pipe. Hence, for the purpose of {optimal experiment design}, these groups are discarded.  As a result, we consider $n_{\theta}=11$ groups in \texttt{BWFLnet}, as showed in Figure~\ref{fig:bwflnet_groups}.



Pressure sensors are installed at inlets and outlets of PCVs and DBVs, the top water inlet, and at the two critical points, defined as the nodes with the lowest pressure. This results in $13$ pressure sensors whose locations are fixed. Other $24$ pressure sensors have been placed according to expert engineering judgment - see Figure~\ref{fig:bwflnet}. In this study, we solve {optimal experiment design} problems in \texttt{BWFLnet}, and compare the obtained solutions with the expert choice. Within the problem formulation, we ignore the expert choice, and we incorporate the $13$ fixed sensor locations by setting to $1$ the lower bounds of the corresponding entries in $z \in \mathbb{R}^{n_n}$. In addition, matrix $X_0$ is defined by the sensitivity vectors associated with the installed flow sensors at PCVs, DBVs and water inlets. Analogously to previous literature (e.g. \cite{Kapelan2005}), we formulate the problem of {optimal experiment design} considering three different demand conditions ($n_t=3$): morning peak ($07:00$), afternoon dip ($14:15$), and evening peak ($18:00$). An alternative strategy could utilize demand levels corresponding to fire flow tests, where fire hydrants are opened to increase flow velocities across network pipes \cite{Walski1983}. However, many water utilities (e.g. in the UK) do not conduct fire flow tests as these tests are resource intensive, time-consuming and potentially disruptive (e.g. discoloration events) - see the discussion in \cite{Waldron2021Disc}. Hydraulic equations~\eqref{eq:hyd_eqs} are solved using the null-space solver developed by \cite{Abraham2015}. All sensitivity vectors $\nabla \phi_{j,k}$ are computed by the computationally efficient method outlined in supplementary material S1.
\subsection{Single-objective optimization}
We first implement the developed convex heuristic to individually minimize $f_D$ and $f_T$ in \texttt{BWFLnet}, i.e. solve Problem~\eqref{eq:prob_indiv}. The number of sensors to be added to the fixed $13$ measuring devices ranges from $2$ to $40$. Hence, we consider the installation of $15$ to $53$ pressure sensors. We formulate Problem~\eqref{eq:prob_indiv} for each value of $m \in \{15,\ldots,53\}$, and implement the developed convex heuristic. As a first step, we solve a convex relaxation~\eqref{eq:prob_convex} of the resulting optimization problem. Then, we implement the Round-and-Swap algorithm to compute a good quality integer solution. The implemented method returns both lower and upper bounds to the optimal value of the original mixed-integer optimization problem. 

\begin{figure}
    \centering
    \subfloat[Upper and lower bounds. \label{fig:bwflnet_dopt_obj}]{\includegraphics[width=0.48\textwidth]{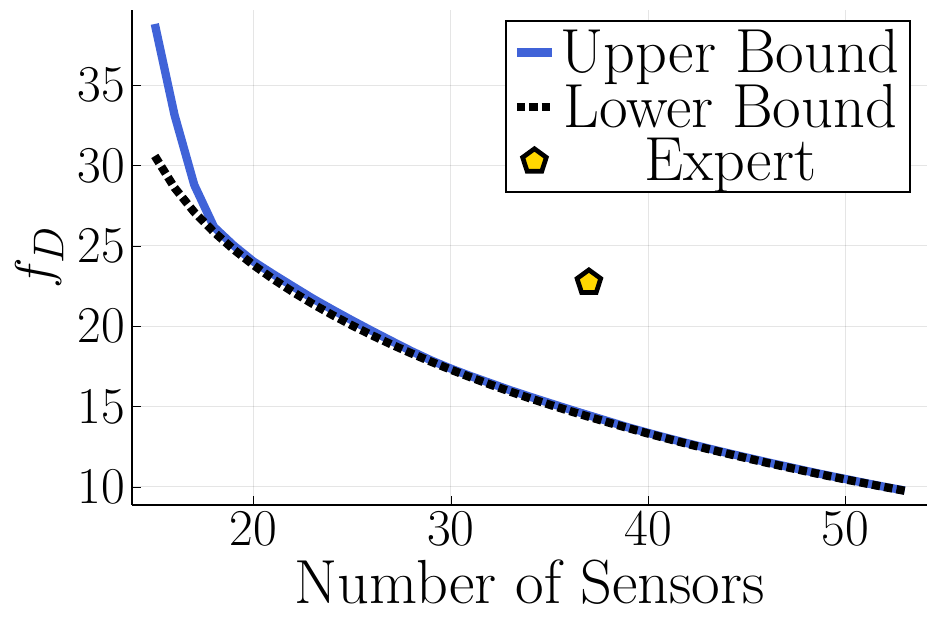}}
    \,
    \subfloat[Bounds on optimality gaps. \label{fig:bwflnet_dopt_gap}]{\includegraphics[width=0.48\textwidth]{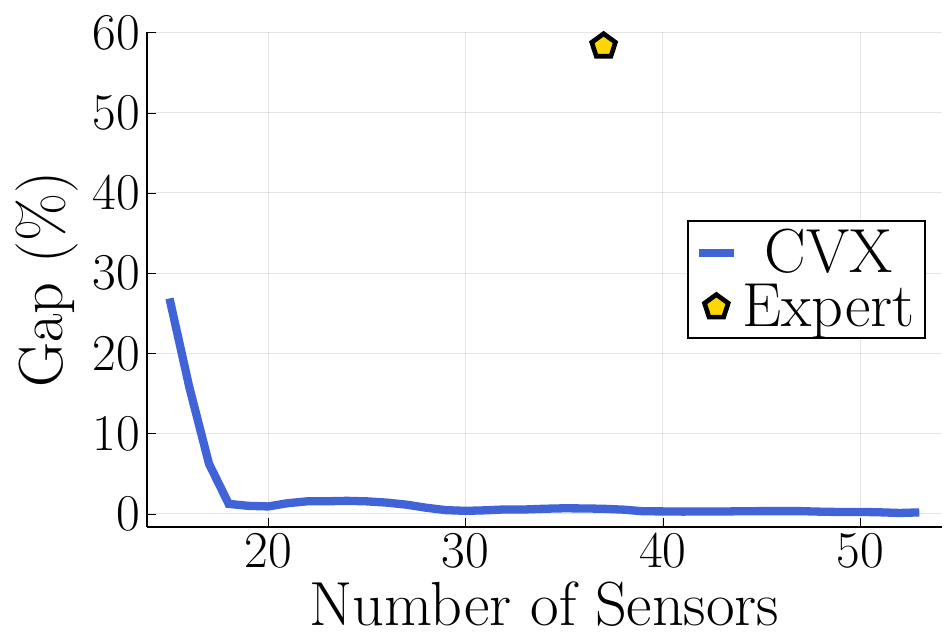}}
    \caption{Results obtained by the convex heuristic (CVX) when implemented to install $15$ to $53$ pressure sensors in \texttt{BWFLnet} for minimizing $f_D$.}
    \label{fig:bwflnet_dopt}
\end{figure}

\begin{figure}
    \centering
    \subfloat[Upper and lower bounds. \label{fig:bwflnet_pmedian_obj}]{\includegraphics[width=0.48\textwidth]{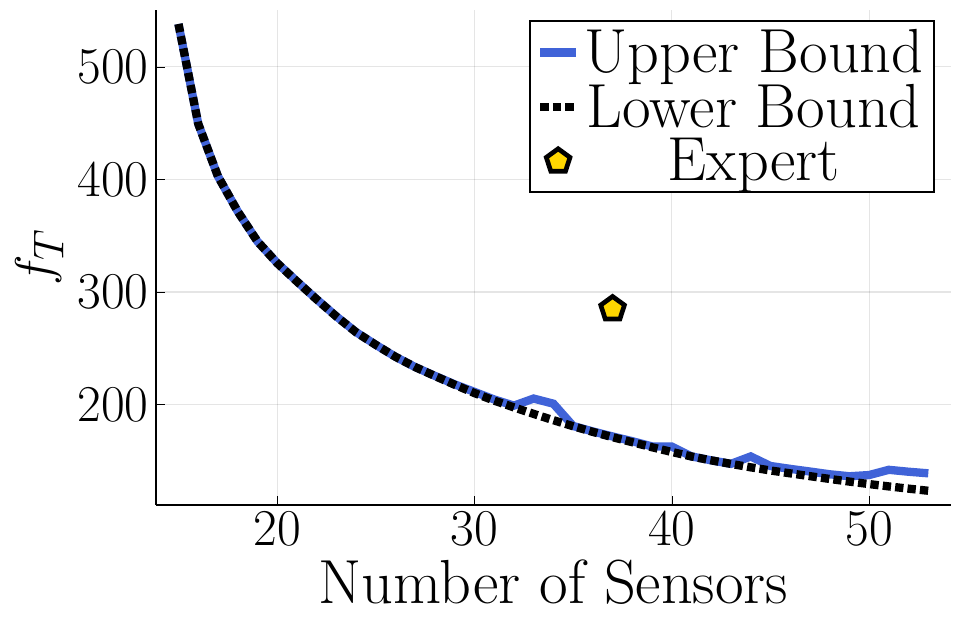}}
    \,
    \subfloat[Bounds on optimality gaps. \label{fig:bwflnet_pmedian_gap}]{\includegraphics[width=0.48\textwidth]{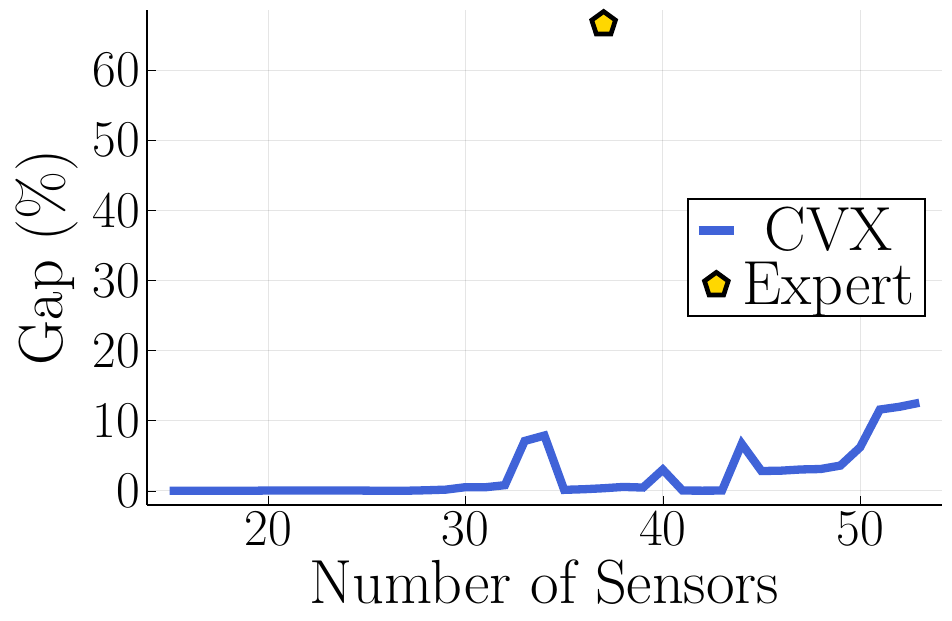}}
    \caption{Results obtained by the convex heuristic (CVX) when implemented to install $15$ to $53$ pressure sensors in \texttt{BWFLnet} for minimizing $f_T$.}
    \label{fig:bwflnet_pmedian}
\end{figure}

The results obtained for $f_D$ are reported in Figure~\ref{fig:bwflnet_dopt}, and we observe that the gap between lower and upper bounds is small for most instances. Therefore, with the exceptions of problem instances with two and three additional sensors, all computed feasible solutions are proven to be globally optimal or near-optimal. While the results for other problem instances suggest that the solutions for $m \in \{15,16\}$ are near-optimal, the lower bounds for these instances are not as tight.
Similarly, in the case of $f_T$, the implemented convex heuristic successfully computes near-optimal feasible solutions for the majority of problem instances - see Figure~\ref{fig:bwflnet_pmedian}. Figure~\ref{fig:bwflnet_pmedian_gap} shows that computed lower and upper bounds have larger differences only for few problem instances, suggesting that the Round-and-Swap algorithm has converged to sub-optimal solutions in these cases. 

Figures \ref{fig:bwflnet_dopt} and \ref{fig:bwflnet_pmedian} 
enable a cost-benefit analysis of different instrumentation levels, as well as compare the existing sensor configuration against optimal solutions. In particular, these results show that the computed solutions with the same number of sensors as the expert choice (i.e. $m=37$) dominate the existing sensor configuration. In fact, in the case $f_D$, the convex heuristic results in a bound on the optimality gap equal to $0.64\%$, while the bound on the gap of the expert solution is $58\%$. In the case of $f_T$, the heuristic has converged to a solution that is within $0.36\%$ of optimality, while the expert solution has a bound on the optimality gap equal to $66\%$. 


Although the considered sensor placement problem is solved off-line,
and there are no constraints on computational speed and time,
the analysis of the computational effort provides insights into the
implemented problem formulations and algorithms. In the case of D-optimality, the convex Problem~\eqref{eq:prob_convex} is rewritten as a semidefinite program - see Problem~\eqref{eq:dopt_convex} in \ref{sec:reform}. Even though the problem includes thousands of variables, the convex solver MOSEK has successfully converged in all tested instances with limited computational effort - the maximum computational time was $0.55$ seconds. Analogously, the Round-and-Swap algorithm converged in less than {$300$} iterations (i.e. function evaluations) for all problems. As a result, the longest computational time was $0.07$ seconds. In the case of $f_T$, we report a more challenging computational performance. In fact, the formulation of Problem~\eqref{eq:pmedian_convex} results in a linear program with $4,935,062$ variables. The state-of-the-art solver GUROBI converged in all problem instances, with a maximum computational time of $250$ seconds. In comparison, the Round-and-Swap algorithm took less { than} $500$ iterations in the majority of tested instances, but some problem formulations required up to $15,164$ iterations (i.e. function evaluations). As a consequence, the maximum computational time was $190$ seconds.

\begin{figure}
    \centering
    \subfloat[Minimization of $f_D$. \label{fig:bwflnet_sensors_dopt}]{\includegraphics[width=0.48\textwidth]{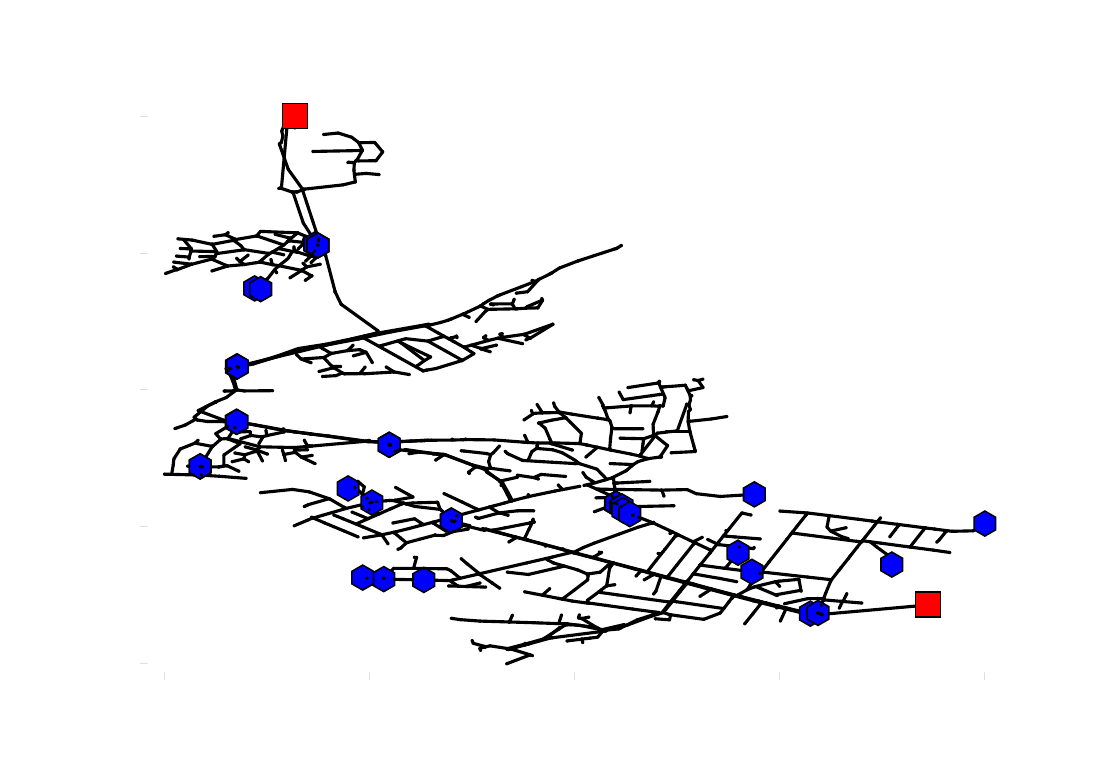}}
    \,
    \subfloat[Minimization of $f_T$. \label{fig:bwflnet_sensors_pmedian}]{\includegraphics[width=0.48\textwidth]{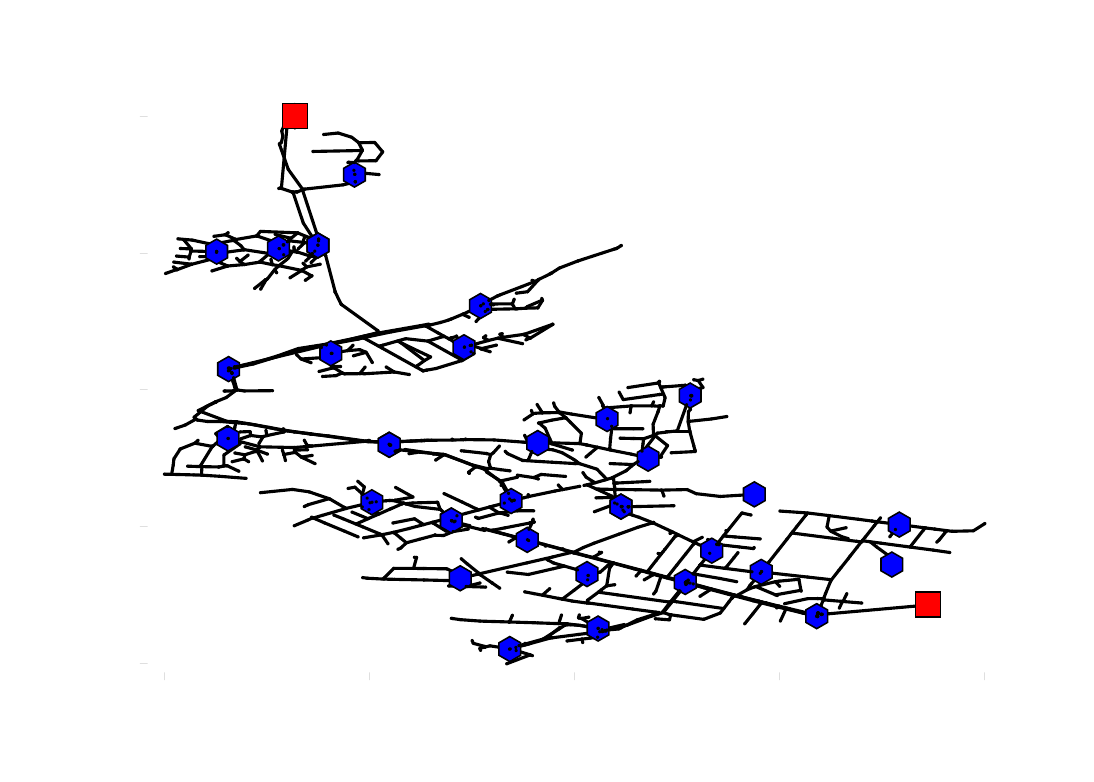}}
    \caption{Optimized valve configurations for individual objectives with $m=37$.}
    \label{fig:bwflnet_sensors_indiv}
\end{figure}

In {Figures~\ref{fig:bwflnet_sensors_dopt} and \ref{fig:bwflnet_sensors_pmedian}}, we report the sensor configurations obtained when minimizing $f_D$ and $f_T$, respectively, when the number of sensors to install is $m=37$. As expected, when $f_D$ is minimized, we do not obtain an even spread of sensors, which appear to be grouped around certain network links to minimize the uncertainty around the estimated parameters. However, pressure measuring devices are commonly used for additional tasks, for example leak localization~\cite{Blocher2021}. In such cases, a more uniform sensor distribution is preferred. This is obtained by minimizing $f_T$, as shown in Figure~\ref{fig:bwflnet_sensors_pmedian}.

\subsection{Bi-objective optimization}
\begin{figure}
    \centering
    \includegraphics[width=0.7\textwidth]{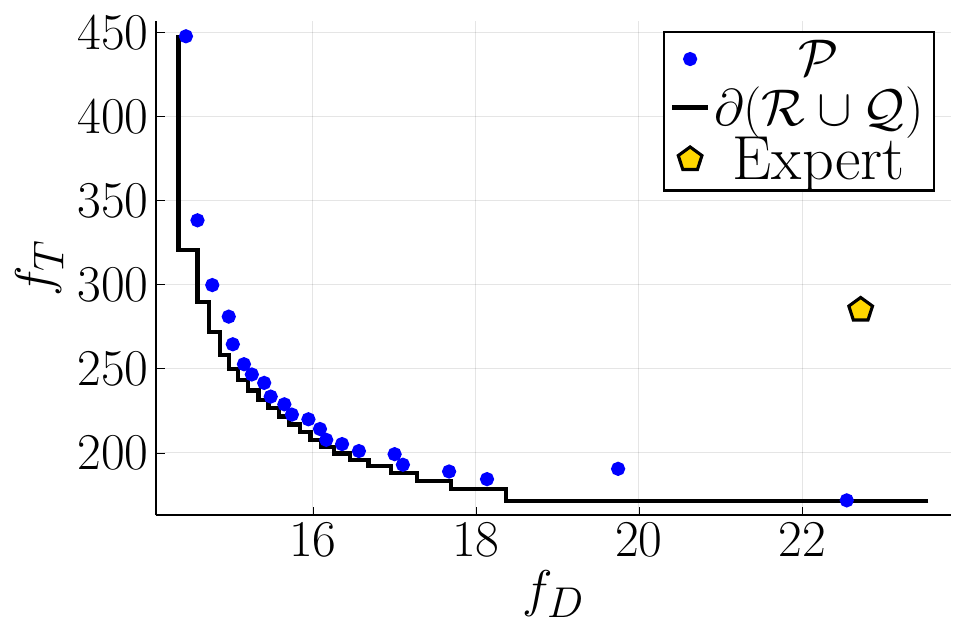}
    \caption{Approximation of the Pareto front obtained by Algorithm~\ref{alg:cheb} with $N=20$, implemented to solve the bi-objective sensor location problem of installing $37$ sensors in \texttt{BWFLnet}.}
    \label{fig:bwflnet_front}
\end{figure}

An ideal sensor configuration corresponds to a compromise between accuracy of hydraulic model calibration, and uniform spatial distribution. In order to investigate sensor configurations representing trade-offs between $f_D$ and $f_T$, we formulate the bi-objective Problem~\eqref{eq:prob_multiobj} to optimize the locations of $m=37$ measuring devices, including the $13$ fixed sensors. We implement the Chebyshev scalarization scheme in Algorithm~\ref{alg:cheb} with $N=20$, where single-objective optimization problems are solved by { the convex heuristic} to compute feasible solutions and lower bounds on the optimal values - see Section~\ref{sec:singleopt}. The algorithm generates a set of potentially non-dominated points $\mathcal{P}$, which corresponds to $22$ feasible solutions. This set is an approximation of the true Pareto front, computed by the convex heuristic. In addition, Algorithm~\ref{alg:cheb} generates sets $\mathcal{R}$ and $\mathcal{Q}$ such that no point of the true Pareto front is included in $\mathcal{R} \cup \mathcal{Q}$ - see Theorem~\ref{thm:bound}. As a result, $\partial(\mathcal{R} \cup \mathcal{Q})$ is a bound on the Pareto front, as it limits the space that does not contain non-dominated points (i.e. points on the Pareto front). This is presented in Figure~\ref{fig:bwflnet_front}, which reports the results computed by Algorithm~\ref{alg:cheb}. The scalarization method has generated an even spread of potentially non-dominated points, which represents the trade-offs between $f_D$ and $f_T$. Algorithm~\ref{alg:cheb} requires the implementation of the convex heuristic to solve $20$ instances of Problem~\eqref{eq:prob_cheb}, resulting in a total computational time of roughly $4.5$ hours. The computational effort is dominated by the solution of the convex relaxation, which is rewritten as Problem~\eqref{eq:cheb_convex3}. This results in a semidefinite program with $4,935,073$ scalar variables, one matrix variable and $11$ exponential cone constraints. The average computational time for MOSEK to converge to a solution was $812$ seconds. In comparison, the largest number of iterations in the Round-and-Swap algorithm was $1522$, which resulted in a maximum computational time of approximately $16$ seconds.

\begin{figure}
    \centering
    \subfloat[Pareto front with highlighted solution. \label{fig:bwflnet_tradeoff_front}]{\includegraphics[width=0.48\textwidth]{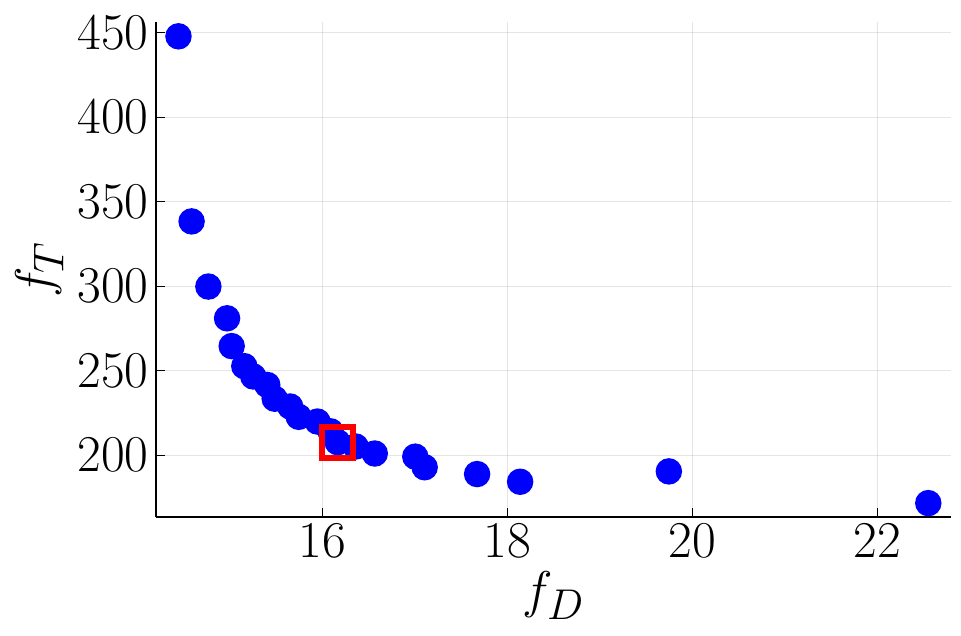}}
    \,
    \subfloat[Corresponding sensor configuration. \label{fig:bwflnet_tradeoff_sensors}]{\includegraphics[width=0.48\textwidth]{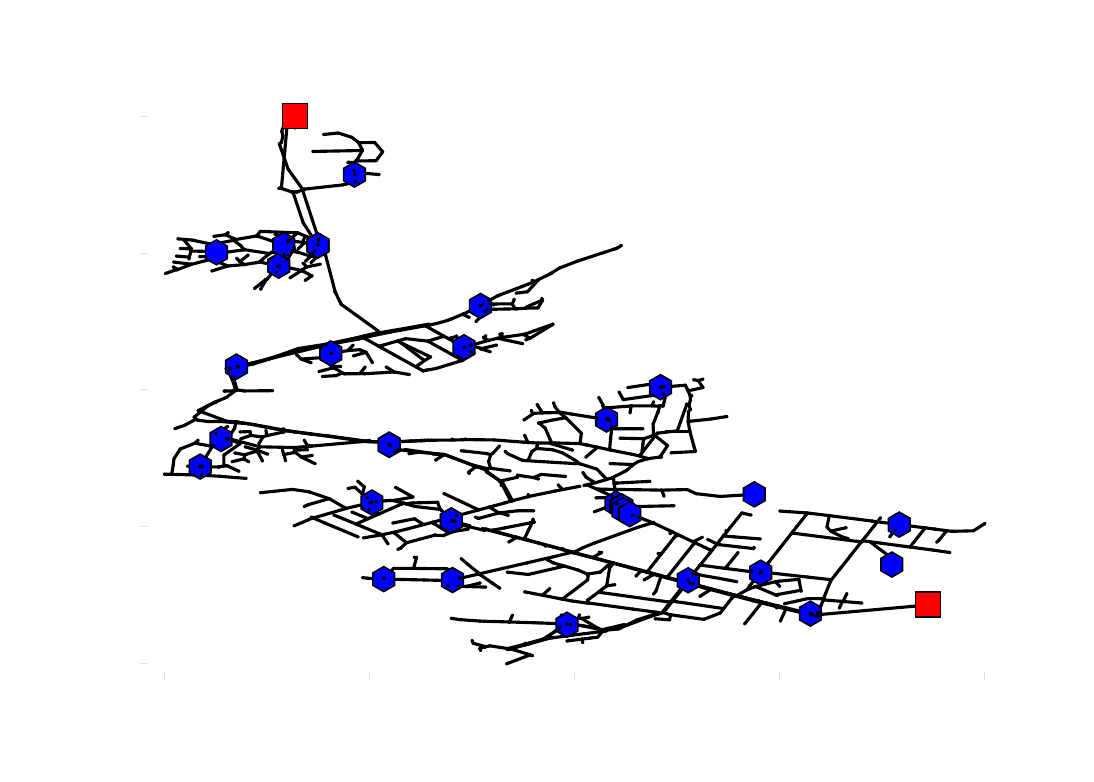}}
    \caption{Point on the Pareto front with corresponding sensor configuration.}
    \label{fig:bwflnet_tradeoff}
\end{figure}

We observe that, for the considered case study, sets $\mathcal{W}$ and $\mathcal{R} \cup \mathcal{Q}$ in Theorems~\ref{thm:upbound} and \ref{thm:bound}, respectively, results in tight bounds on the Pareto front. This suggests that set $\mathcal{P}$ computed by Algorithm~\ref{alg:cheb} is close to the true Pareto front of Problem~\eqref{eq:prob_multiobj}. Figure~\ref{fig:bwflnet_tradeoff} illustrates the trade-offs between $f_D$ and $f_T$. The highlighted point on the approximated Pareto front in Figure~\ref{fig:bwflnet_tradeoff_front} represents a compromise between sensitivity and topology based optimality definitions. As shown in Figure~\ref{fig:bwflnet_tradeoff_sensors}, this results in a sensor configuration with a better spatial coverage compared to Figure~\ref{fig:bwflnet_sensors_dopt}. Therefore, such a solution meets other operational objectives in addition to hydraulic model calibration, including detection and localization of water leaks.

Finally, we have tested the proposed methods on other network models and no significant difference was observed. Supplementary material S4 reports the results obtained for a different case study, the water network model \texttt{L-TOWN-A}~\cite{Vrachimis2020}. {Similarly to what occurs for the case study \texttt{BWFLnet}}, Algorithm~\ref{alg:cheb} is able to approximate the Pareto fronts with tight bounds. However, we decided to not include these results in the main text as they do not provide any additional insight.

\subsection{Numerical experiments with IPOPT}
\label{app:ipopt}
We have investigated the implementation of the state-of-the-art non-linear programming solver IPOPT~\cite{Wachter2006} to solve the considered convex problems without the need for semidefinite programming reformulations, i.e. directly considering Problems~\eqref{eq:dopt_convex_0} and \eqref{eq:cheb_convex}. We consider the same case study \texttt{BWFLnet} as in Section~\ref{sec:numexp}. In the implementation of IPOPT, we have passed exact Jacobian and Hessian matrices through the low-level C-interface provided by Ipopt.jl (\url{https://github.com/jump-dev/Ipopt.jl}). In the case of the individual minimization of $f_D$ (i.e. Problem~\eqref{eq:dopt_convex_0}), IPOPT required an average computational time of $238$ seconds, compared to the average computational time associated with MOSEK, which was $0.35$ seconds. In the case of $f_{\beta}$ (i.e. Problem~\eqref{eq:cheb_convex}), the implementation of IPOPT resulted in an out-of-memory error. This {poor} computational performance can be explained by the structure of Problems~\eqref{eq:dopt_convex_0} and \eqref{eq:cheb_convex}, which result in large and dense Jacobian and Hessian matrices. As example, in the case of $f_{\beta}$, the number of non-zero elements within the Jacobian matrix is equal to $19,704,782$. These experiments {support} the choice of reformulating these optimization problems as SDPs and solving them using MOSEK. Not only such reformulations do not require the analytical derivation of Jacobian and Hessian matrices, but they also result in faster solution times. 
We also expect off-the-shelf convex mixed-integer non-linear programming solvers like BONMIN \cite{Bonami2008} to have a {poor} computational performance when applied to the considered MIPs, because they rely on IPOPT as non-linear programming solver. 
\subsection{Application of MOEAs}
\label{sec:nsgaii}
\begin{figure}[h!]
    \centering
    \subfloat[Complete Pareto front approximations.]{\includegraphics[width=0.48\textwidth]{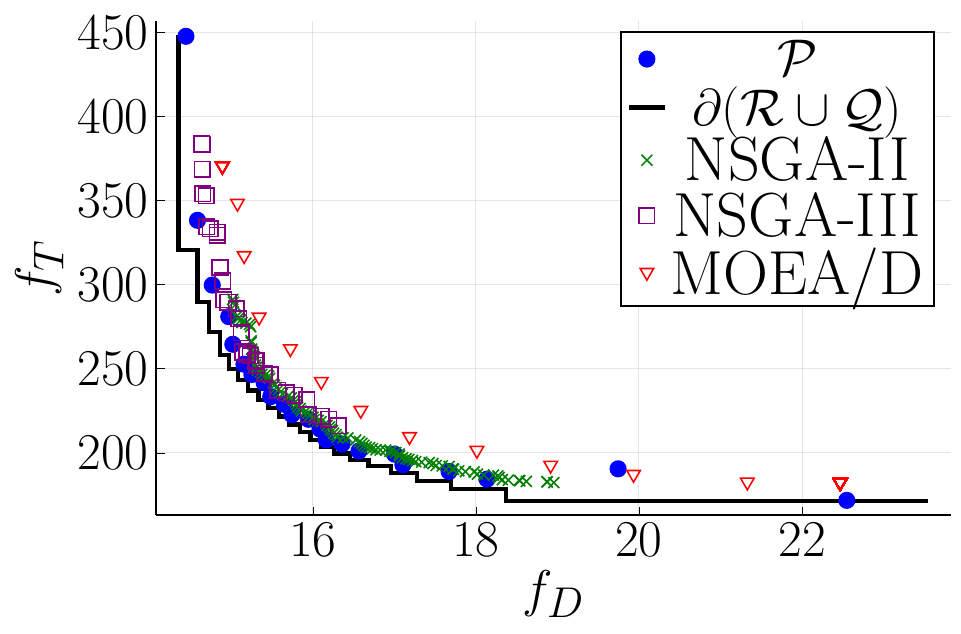}}
    \,
    \subfloat[Zoomed image.]{\includegraphics[width=0.48\textwidth]{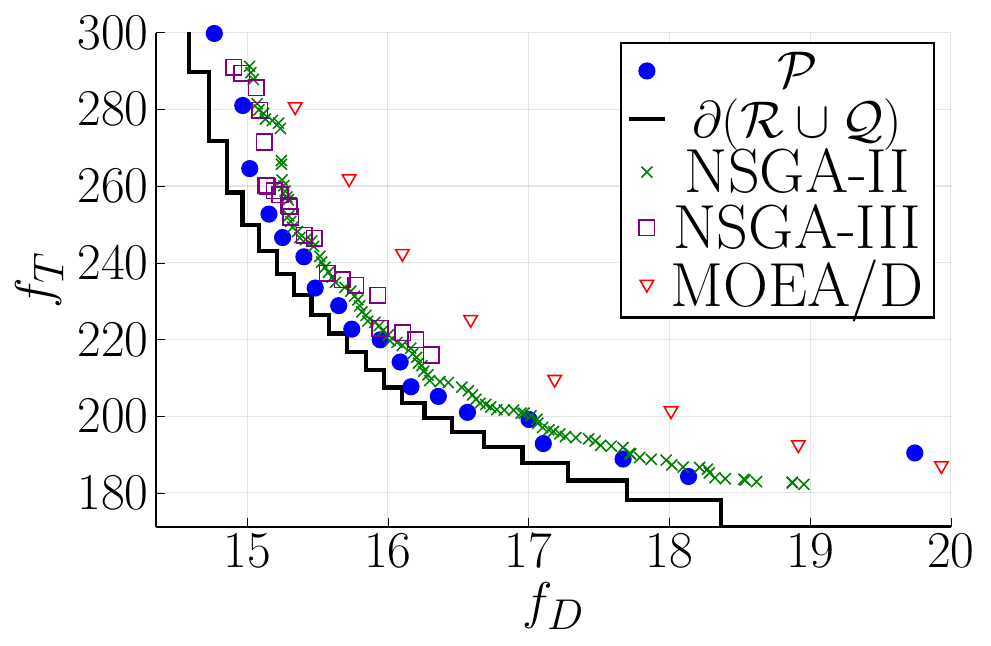}}
    \caption{Comparison of Algorithm~\ref{alg:cheb} and NSGA-II, NSGA-III, and MOEA/D.}
    \label{fig:bwflnet_moeas}
\end{figure}

{We compare the Pareto front approximation computed by Algorithm~\ref{alg:cheb} with those obtained by off-the-shelf Multi-objective Optimization Evolutionary Algorithms (MOEAs). The considered solvers are NSGA-II~\cite{Deb2002}, NSGA-III~\cite{Deb2013}, and MOEA/D~\cite{Zhang2007}}. In the implementation of MOEA/D, we select Chebyshev scalarization functions, analogously to Algorithm~\ref{alg:cheb}. The solvers are implemented through Platypus (\url{https://github.com/Project-Platypus/Platypus}). Because Algorithm~\ref{alg:cheb} required a computational time of $4.5$ hours, we set a time limit of $6$ hours. As shown in Figure~\ref{fig:bwflnet_moeas}, the solutions computed by the MOEAs are dominated by those obtained with the convex heuristic. Moreover, Algorithm~\ref{alg:cheb} results in additional information compared to MOEAs. In fact, not only it computes a set of potentially non-dominated solutions, but also theoretical bounds on the true Pareto front.

\section{Conclusions}
We have proposed a new bi-objective problem formulation for optimal placement of pressure sensors in water networks. The considered objectives are the minimization of uncertainty in hydraulic model parameters and the maximization of spatial coverage. 
The resulting bi-objective optimization problem belongs to the class of convex mixed-integer non-linear programming. We have considered a Chebyshev scalarization method, which approximates the Pareto front solving a sequence of convex mixed integer non-linear programs. Computing exact solutions for these optimization problems can be impractical when the number of integer variables is large. Hence, we have implemented a convex heuristic to generate feasible solutions with lower bounds on the globally optimal values. Furthermore, we have demonstrated that such lower bounds can be used to define a region containing the true Pareto front of the bi-objective optimization problem. These results are not restricted to the framework of {optimal experiment design} and they can inform researchers working on bi-objective optimization in other application areas. 
We evaluate the proposed methods using an operational water network from the UK with over $2000$ nodes. The developed convex heuristic computes near-optimal solutions for the individual objective minimization problems, and tight bounds on the true Pareto front of the bi-objective optimization problem. 
The overall computational time required to generate $20$ points on the (approximated) Pareto front was less than $5$ hours, the vast majority of which was spent solving convex relaxations. This suggests that the convex heuristic can be used to solve bi-objective sensor placement problems arising in practical applications. Hence, the presented methodology is a stepping stone in the development of computational tools to support the upgrade of existing infrastructure with the implementation of advanced monitoring and control solutions. Future work should extend this framework to take into account the different sources of uncertainty inherent in hydraulic models, including the status of partially closed valves, pipes' lengths and diameters, and estimated customer demands. 

\appendix
\section{Efficient computation of sensitivity vectors}
\label{sec:senscomp}
We develop a computationally efficient method to calculate the vectors of sensitivities $\nabla \phi_{ik}(\bar{\theta})$ for all $i \in \{1,\ldots,n_p+n_n\}$, $k \in \{1,\ldots,n_t\}$. We have that
\begin{equation}
     \nabla \phi_{ik}(\bar{\theta}) = J_k^Te_i,
\end{equation}
where $J_k$ is the Jacobian matrix of function $\phi_k$ evaluated at $\bar{\theta}$, and $e_i$ is a vector of the canonical basis of $\mathbb{R}^{n_p+n_n}$. Recall that $\phi_k(\bar{\theta})=[\bar{q}_k^T\,\bar{h}_k^T]^T$ is such that
\begin{equation}
    \label{eq:hyd_eqs}
    \begin{split}
        &A_{11}(\bar{q}_k,\bar{\theta})\bar{q}_k + A_{12} \bar{h}_k + A_{10}{h_0}_k + A_{13}\eta_k = 0\\
        &A_{12}^T\bar{q}_k = d_k
    \end{split}
\end{equation}
Hence, by the implicit function theorem, we have that the Jacobian of function $\phi_k$ is given by:
\begin{equation}
\label{eq:implicit}
    J_k = - \begin{bmatrix}
    NA^k_{11} & A_{12} \\
    A_{12}^T & 0
    \end{bmatrix}^{-1}\begin{bmatrix}B^k\\0\end{bmatrix},
\end{equation}
where $N=\textbf{diag}(n_1,\ldots,n_{n_p})$, $A_{11}^k=A_{11}(\bar{q}_k,\bar{\theta})$, and $B^k \in \mathbb{R}^{n_p\times n_r}$ is defined by: \begin{equation}
    B^k_{lj} = \frac{\partial \rho_l(\bar{\theta})}{\partial \theta_j}|\bar{q}_{lk}|^{n_l-1}\bar{q}_{lk}, \quad  l \in \{1,\ldots,n_p\}, \; j \in \{1,\ldots,n_r\}.
\end{equation}
Hence, we have:
\begin{equation}
    \nabla \phi_{ik}(\bar{\theta})  = -\begin{bmatrix} (B^k)^T & 0 \end{bmatrix}\begin{bmatrix} NA_{11}^k & A_{12} \\ A_{12}^T & 0 \end{bmatrix}^{-1}e_i
\end{equation}
The sensitivity vector is given by:
\begin{equation}
   \nabla \phi_{ik}(\bar{\theta}) = -(B^k)^T \delta_q^{(ik)},
\end{equation}
with vectors $\delta_q^{(ik)} \in \mathbb{R}^{n_p}$ and $\delta_h^{(ik)}\in \mathbb{R}^{n_n}$ solutions of the linear system:
\begin{equation}
\label{eq:saddle}
    \begin{bmatrix} NA_{11}^k & A_{12} \\ A_{12}^T & 0 \end{bmatrix}\begin{bmatrix}\delta_q^{(i,k)}\\ \delta_h^{(i,k)}\end{bmatrix} = \begin{bmatrix}e_{1i} \\ e_{2i} \end{bmatrix}
\end{equation}
where $e_i = [e_{1i}^T \, e_{2i}^T]^T$. Hence, to compute the sensitivity vectors, we need to solve $n_t(n_p+n_n)$ linear systems like \eqref{eq:saddle}. This would require computing $n_t$ different factorizations for the matrix in \eqref{eq:saddle}, which has order $n_p+n_n$. When large water networks are considered (i.e. $n_p>1000$), this results in significant computational effort. However, system \eqref{eq:saddle} has a saddle point structure \cite{Benzi2005}. Similar matrices are also encountered when solving \eqref{eq:hyd_eqs} using the Newton-Raphson method. Previous work has proposed a null-space algorithm to solve saddle point systems analogous to~\eqref{eq:saddle} at each Newton-Raphson iteration - see \cite{Abraham2015}. Here, we take a similar approach, and we derive a null-space method to solve \eqref{eq:saddle}. 
 
The incidence matrix $A_{12}$ has full column rank in the case of water distribution networks \cite{Elhay2014}. Moreover, we assume $\text{Ker}(A_{11}^k)\cap \text{Ker}(A_{12}^T)=\{0\}$, which is equivalent to requiring that the water network does not include loops where all flows are equal to zeros. Let $Z \in \mathbb{R}^{n_p \times (n_p - n_n)}$ be null-space basis of $A_{12}^T$. Then, any vector $\delta_q^{(ik)}$ satisfying~\eqref{eq:saddle} is written as $\delta_q^{(ik)} = A_{12}w + Zv$, for $w \in \mathbb{R}^{n_n}$ such that $A_{12}^TA_{12}w = e_{2i}$ and $v \in \mathbb{R}^{n_p-n_n}$. By substituting for $\delta_q^{(ik)}$ into the first block of equations in \eqref{eq:saddle}, and pre-multiplying by $Z^T$, we obtain:
\begin{equation}
\label{eq:nullspace_system}
    Z^TNA_{11}^kZv = Z^T(e_{1i}-NA_{11}^kA_{12}w),
\end{equation}
where we used the fact that $Z^TA_{12}=0$. Instead of the linear system \eqref{eq:saddle}, we now have to solve the smaller system of linear equations \eqref{eq:nullspace_system}. The proposed method requires computing factorizations of $A_{12}^TA_{12}$ and $Z^TNA_{11}^kZ$, which are symmetric positive definite matrices (SPD) of orders $n_n$ and $n_p-n_n$, respectively. In fact, since $A_{12}$ has full rank, $A_{12}^TA_{12}$ is SPD. Moreover, $\text{Ker}(A_{11}^k)\cap \text{Ker}(A_{12}^T)=\{0\}$ implies that $Z^TNA_{11}^kZ$ is SPD. The overall algorithm to compute sensitivity vectors $\nabla \phi_{ik}(\bar{\theta})$ for all $i \in \{1,\ldots,n_p+n_n\}$ and $k \in \{1,\ldots,n_t\}$ is outlined in Algorithm \ref{alg:senscomp}. It is necessary to compute a single Cholesky factorization of $A_{12}^TA_{12}$, and $n_t$ different Cholesky factorizations for $Z^TNA_{11}^kZ$. All linear systems are then solved by forward and back substitutions.
\begin{algorithm}
\begin{algorithmic}
\State{\textbf{Input:} matrices $A_{12},A_{11}^1,\ldots,A_{11}^{n_t},N$}
\State{\textbf{Output:} $\nabla \phi_{ik}(\bar{\theta})$  for all $i \in \{1,\ldots,n_p+n_n\}$, $k \in \{1,\ldots,n_t\}$.}
\State{Compute a null-space basis of $A_{12}^T$ denoted by $Z \in \mathbb{R}^{n_p \times (n_p-n_n)}$.}
\State{Compute a Cholesky factorization $A_{12}^TA_{12}=LL^T$.}
\For{$i \in \{1,\ldots,n_p+n_n\}$}
\State{Let $e_i=[e_{1i}^T \, e_{2i}^T]^T$ be the $i$-th vector of the canonical basis of $\mathbb{R}^{n_p+n_n}$.}
\State{Let $w_i \in \mathbb{R}^{n_n}$ be the solution of $LL^Tw_i = e_{2i}$.}
\EndFor
\For{$k=1,\ldots,n_t$}
\State{Compute a Cholesky factorization $Z^TNA_{11}^kZ=P^k(P^k)^T$.}
\For{$i =1,\ldots,n_p+n_n$}
\State{Solve $P^k(P^k)^Tv^k_i=Z^T(e_{1i}-NA_{11}^kA_{12}w_i)$.}
\State{Set $\delta_q^{(ik)} = A_{12}w_i + Zv^k_i$.}
\State{$\nabla \phi_{ik}(\bar{\theta})  = -(B^k)^T\delta_q^{(ik)}$}
\EndFor
\EndFor
\end{algorithmic}
\caption{Null-space method for computing sensitivity vectors}
\label{alg:senscomp}
\end{algorithm}

\cleardoublepage
\section{Technical proofs}
\subsection*{Proof of Lemma 2.1}
    To show this, let $u \in \mathbb{R}^{n_n\times n_n}$ be a feasible for Problem~(9) in the main manuscript. We have:
    \begin{equation}
        \begin{split}
        \sum_{i=1}^{n_n}\sum_{j = 1}^{n_n}C_{ij}u_{ij} &=\sum_{i =1}^{n_n}\sum_{j \in N_i \cap I(\hat{z})}C_{ij}u_{ij}\\
        &\geq \sum_{i =1}^{n_n}\min_{j \in N_i \cap I(\hat{z})}C_{ij}\Bigg(\sum_{j \in N_i \cap I(\hat{z})}u_{ij}\Bigg)\\
        &=\sum_{i =1}^{n_n}\min_{j \in N_i \cap I(\hat{z})}C_{ij}\Bigg(\sum_{j =1}^{n_n}u_{ij}\Bigg)\\
        &=f_T(\hat{z})
        \end{split}
    \end{equation}
    Let $j^* \in \text{argmin}\big\{C_{ij}\, | \, j \in N_i \cap I(\hat{z})\big\}$. Then, define $\hat{u} \in \mathbb{R}^{n_n \times n_n}$ by setting:
    \begin{equation}
        \hat{u}_{ij} = \begin{cases}
            1 & \text{if } j = j^*\\
            0 & \text{otherwise}.
        \end{cases}, \quad \forall j \in \{1,\ldots,n_n\},
    \end{equation}
    for each $i \in \{1,\ldots,n_n\}$. By construction, $\hat{u}$ is feasible for Problem~(9), and
    \begin{equation}
        \begin{split}
        \sum_{i=1}^{n_n}\sum_{j=1}^{n_n}C_{ij}\hat{u}_{ij} &= \sum_{i=1}^{n_n}C_{ij^*}\\
        & = \sum_{i=1}^{n_n}\min_{j \in N_i \cap I(\hat{z})}C_{ij}\\
        & = f_T(\hat{z})
        \end{split}
    \end{equation}
    Hence, $f_T(\hat{z})$ is the optimal value of Problem~(9).\qed    
\subsection*{Proof of Theorem 3.1}
\label{app:proof_weakly}
    Assume that $\hat{z}_{\beta}$ is not weakly Pareto optimal. In this case, there exists $z'$ feasible for Problem~(10) in the main manuscript such that $f_D(z')<f_D(\hat{z}_{\beta})$ and $f_T(z')< f_T(\hat{z}_{\beta})$. Since $w_D^{\beta}>0$, $w_T^{\beta}>0$, and $f^*_D$ and $f^*_T$ are lower bounds for the individual objective functions, we have:
    \begin{equation}
        \begin{split}
            &w_D^{\beta}(f_D(z')-f_D^*)<w_D^{\beta}(f_D(\hat{z}_{\beta})-f_D^*)\\
            &w_T^{\beta}(f_T(z')-f_T^*)< w_T^{\beta}(f_T(\hat{z}_{\beta})-f_T^*).
        \end{split}
    \end{equation} 
    Hence, $f_{\beta}(z')<f_{\beta}(\hat{z}_{\beta})$. This contradicts the assumption that $\hat{z}_{\beta}$ is a globally optimal solution for Problem~(14). Hence, $\hat{z}_{\beta}$ must be a weakly Pareto optimal solution.\qed
\subsection*{Proof of Theorem 3.2}
Define:
    \begin{equation}
        \beta^* = \frac{(f_D(\hat{z}_D)-f_D^*)(f_T(z^*)-f_T^*)}{(f_D(\hat{z}_D)-f_D^*)(f_T(z^*)-f_T^*) + (f_T(\hat{z}_T)-f_T^*)(f_D(z^*)-f_D^*)}
    \end{equation}
    Then, we have:
    \begin{equation}
        w_D^{\beta^*} = \frac{M}{f_D(z^*)-f_D^*}, \quad w_T^{\beta^*} = \frac{M}{f_T(z^*)-f_T^*},
    \end{equation}
    with 
    \begin{equation}
        M = \frac{(f_D(z^*)-f_D^*)(f_T(z^*)-f_T^*)}{(f_D(\hat{z}_D)-f_D^*)(f_T(z^*)-f_T^*) + (f_T(\hat{z}_T)-f_T^*)(f_D(z^*)-f_D^*)}.
    \end{equation}
    Assume that $z^* \in \Theta$ is not a solution of Problem~(14) for $\beta=\beta^*$. Then, there exists $\hat{z}_{\beta^*}$ such that
    \begin{equation}
        \begin{split}
        f_{\beta^*}(\hat{z}_{\beta^*}) &< f_{\beta^*}(z^*)\\
        &= \max\bigg(w_D^{\beta^*}(f_D(z^*)-f_D^*),w_T^{\beta^*}(f_T(z^*)-f_T^*)\bigg)\\
        &= M
        \end{split}
    \end{equation}
    Therefore, we have:
    \begin{equation}
        \max\bigg(w_D^{\beta^*}(f_D(\hat{z}_{\beta^*})-f_D^*),w_T^{\beta^*}(f_T(\hat{z}_{\beta^*})-f_T^*)\bigg)<M.
    \end{equation}
    Which yields:
    \begin{equation}
        \begin{split}
            & \frac{M}{f_D(z^*)-f_D^*} (f_D(\hat{z}_{\beta^*})-f_D^*) < M \\
            & \frac{M}{f_T(z^*)-f_T^*} (f_T(\hat{z}_{\beta^*})-f_T^*) < M.
        \end{split}
    \end{equation}
    Hence:
    \begin{equation}
        \begin{split}
            &  f_D(\hat{z}_{\beta^*}) < f_D(z^*) \\
            &  f_T(\hat{z}_{\beta^*}) < f_T(z^*).
        \end{split}
    \end{equation}
    The above inequalities contradict the assumption that $z^*$ is a Pareto optimal solution for Problem~(10). Hence, $z^*$ must be an optimal solution of Problem~(14) with $\beta =\beta^*$.\qed

\cleardoublepage
\section{Pseudo-code algorithm description}
\begin{algorithm}[h!]
    \begin{algorithmic}
        \State{\textbf{Input:} $z^*$, $A_1,\ldots,A_{n_n}$, $G=[g_1,\ldots,g_{n_n}]$, $e \in \mathbb{R}^{n_p}$, $m$, $\texttt{MaxIter}$.}
        \State{Let $I=\{i_1,\ldots,i_{n_n}\}$ be such that $z^*_{i_1}\geq\ldots\geq z_{i_{n_n}}^*$.}
        \State{Define $I_r = \texttt{reverse}(I)=\{i_{n_n},\ldots,i_1\}$}
        \State{\textbf{Step 1: Rounding:} initialize $\hat{z}=0_{n_n}$ and $k=1$.}
        \While{$\sum_{i=1}^{n_n}\hat{z}_i<m$}
        \If{$\max(G\hat{z} + g_{i_{k}} - b) \leq 0$}{ $\hat{z}_{i_k}=1$.}
        \EndIf
        \State{$k = k+1$.}
        \EndWhile
        \State{\textbf{Step 2: Swapping:} initialize $\texttt{count}=0$.}
        \While{$\texttt{count}\leq \texttt{MaxIter}$}
            \State{Set \texttt{flag}=1 and $S = \{i \in I_r \, | \, \hat{z}_i=1 \text{ and } 0.1\leq z^*_i\leq 0.9\}$.}
            \State{$nS = \{i \in I \, | \, \hat{z}_i=0 \text{ and } 0.1\leq z^*_i\leq 0.9\}$.}
            \For{$i \in S$}
            \For{$j \in nS$}
            \State{$\texttt{count} = \texttt{count}+1$.}
            \State{$z^{\text{test}}=\hat{z}$, $z^{\text{test}}_i=0$, $z^{\text{test}}_j=1$.}
            \If{$\max(Gz^{\text{test}}-b)\leq 0$ and $f(z^{\text{test}})\leq f(\hat{z})$}
            \State{$\texttt{flag} = 0$ and $\hat{z} = z^{\text{test}}$.}
            \EndIf
            \EndFor
            \EndFor
            \If{\texttt{flag}==1}{ break. \textit{\# convergence reached}}
            \EndIf
        \EndWhile
    \end{algorithmic}
    \caption{Round-and-Swap Algorithm}
    \label{alg:round-swap}
\end{algorithm}
\section{\texttt{L-TOWN-A}}
\label{app:ltown}
\begin{figure}
    \centering
    \subfloat[Network layout. \label{fig:ltown}]{\includegraphics[width=0.48\textwidth]{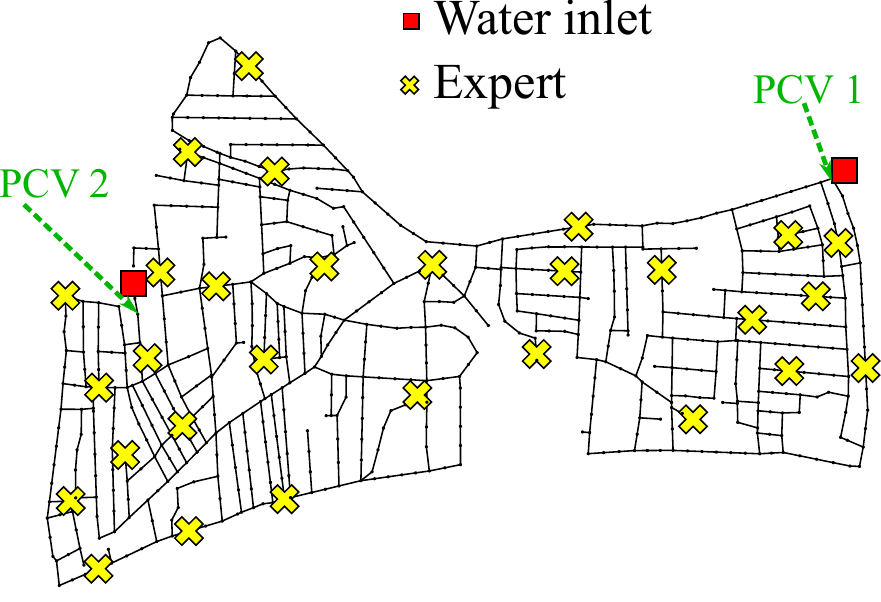}}
    \,
    \subfloat[Pipe groups. \label{fig:ltown_groups}]{\includegraphics[width=0.48\textwidth]{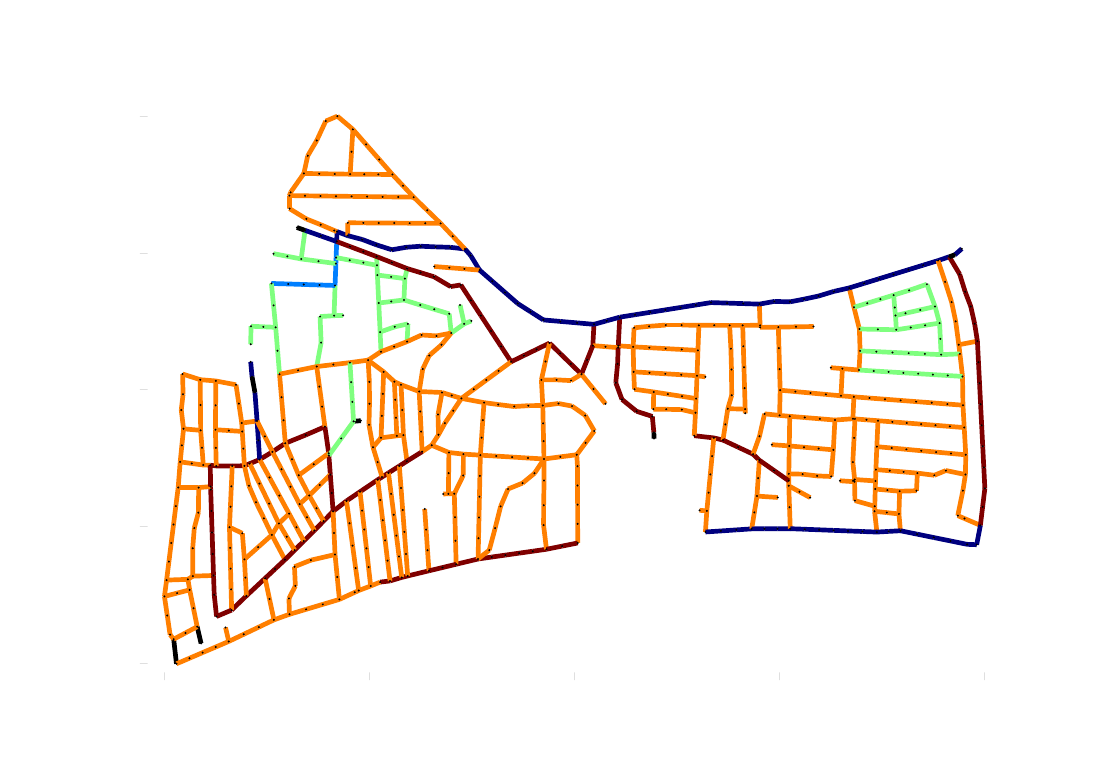}}
    \caption{\texttt{L-TOWN-A}: (a) shows locations of expert choice of existing pressure sensors and pressure control valves (PCVs); (b) pipe group. }
    \label{fig:ltown_all}
\end{figure}

\begin{figure}
    \centering
    \subfloat[Minimization of $f_D$. \label{fig:ltown_dopt}]{\includegraphics[width=0.48\textwidth]{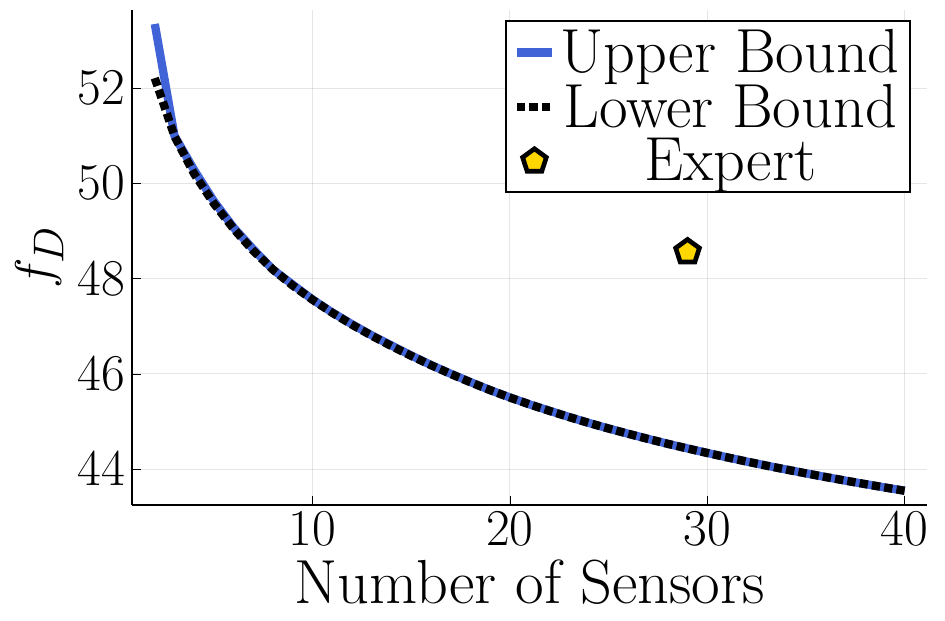}}
    \,
    \subfloat[Minimization of $f_T$. \label{fig:ltown_pmedian}]{\includegraphics[width=0.48\textwidth]{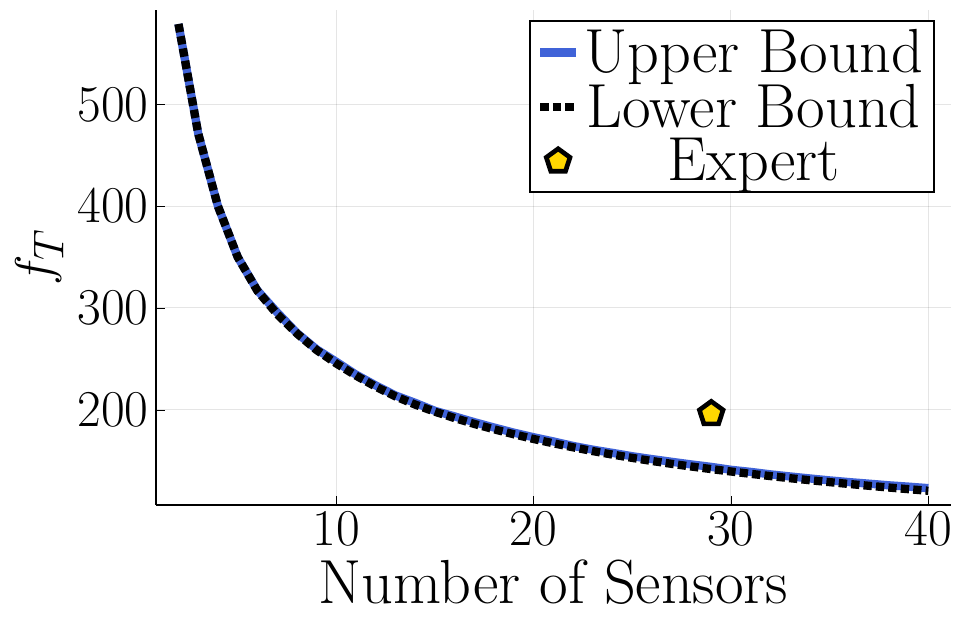}}
    \caption{Results obtained by the convex heuristic implemented to install $2$ to $40$ pressure sensors in \texttt{L-TOWN-A}.}
    \label{fig:ltown_singleopt}
\end{figure}

\begin{figure}[h!]
    \centering
    \includegraphics[width=0.7\textwidth]{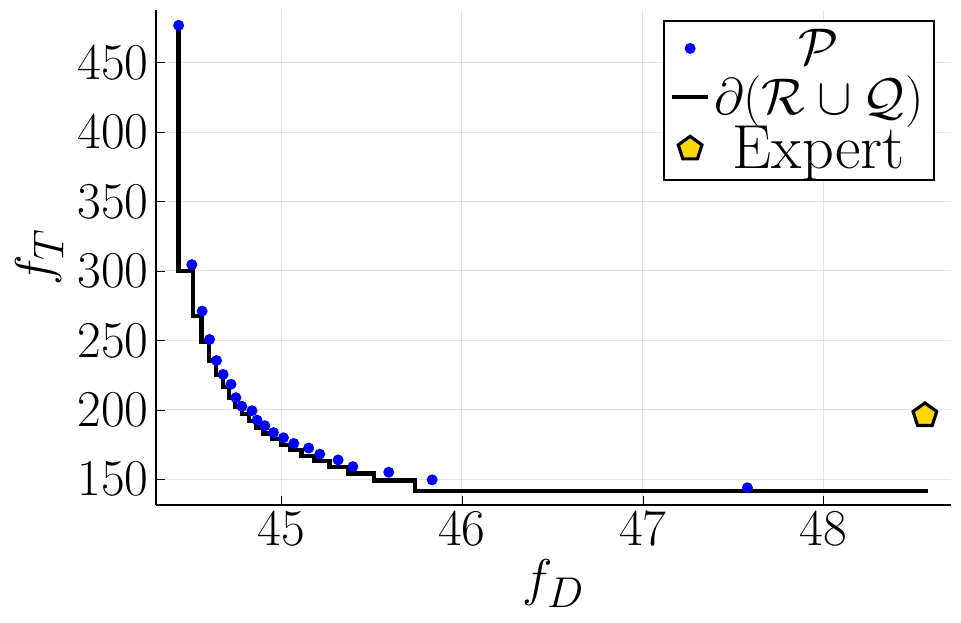}
    \caption{Approximation of the Pareto front obtained by the Chebyshev algorithm with $N=20$, implemented to solve the bi-objective sensor location problem of installing $29$ sensors in \texttt{L-TOWN-A}.}
    \label{fig:ltown_front}
\end{figure}

The water network model L-Town was provided by the organizers of
the BattLeDIM~\cite{Vrachimis2020}, and it represents a hypothetical
town with approximately 10,000 inhabitants, who receive water from two reservoirs. The original network model consists of three distinct hydraulic areas, namely A, B, and C. In this work, we focus on the main network area, referred to as \texttt{L-TOWN-A} and shown in Figure~\ref{fig:ltown}. The remaining smaller areas B and C are fed from \texttt{L-TOWN-A} through a pressure control
valve and a pumping station, respectively. In the model of \texttt{L-TOWN-A}, we opportunely modify the demands at boundary nodes with areas B and C to represent the water use in these sectors. \texttt{L-TOWN-A} includes $768$ links, $661$ nodes, and two pressure control valves installed downstream of the two reservoirs. 
We have grouped pipes in \texttt{L-TOWN-A} in $6$ groups based on their lengths and diameters - the spatial distribution of the groups is shown in Figure~\ref{fig:ltown_groups}. Similarly to \texttt{BWFLnet}, our formulation considers three different demand conditions, corresponding to morning peak, afternoon dip, and evening peak.
The organizers have selected $29$ locations for the installation of pressure sensors in \texttt{L-TOWN-A}, which are presented in Figure~\ref{fig:ltown}. In addition, flow sensors are installed at the two pressure control valves. Pressure and flow sensors have been used in BattLeDIM to detect and localize leaks in \texttt{L-TOWN-A}. 
As firs step, we have formulated and solved two single-objective optimization problems in \texttt{L-TOWN-A} to minimize $f_D$ and $f_T$, respectively. Analogously to \texttt{BWFLnet}, the convex heuristic has resulted in certified globally optimal or near-optimal solutions with tight lower bounds, see Figure~\ref{fig:ltown_singleopt}.
Then, the Chebyshev scalarization algorithm was implemented to approximate the Pareto front of the bi-objective optimization problems to jointly minimize $f_D$ and $f_T$ in \texttt{L-TOWN-A}, with number of sensors to be install equal to $29$. Figure~\ref{fig:ltown_front} shows that the algorithm has successfully generated $20$ potentially non-dominated solutions with a tight bound on the true Pareto front. These solutions dominate the expert choice, and represent the best compromises between sensitivity and topology based metrics.
\section{Reformulation of convex relaxations}
\label{sec:reform}
First, as observed in Lemma~\ref{lem:pmedian_equivalently}, when $f(z)=f_T(z)$, Problem~\eqref{eq:prob_convex} is equivalent to the following linear program with $n_n^2 + n_n$ scalar variables:
\begin{equation}
    \label{eq:pmedian_convex}
    \begin{alignedat}{3}
    &\text{minimize}&\; \; & \sum_{i=1}^{n_n}\sum_{j=1}^{n_n} C_{ij}u_{ij}\\
    &\text{subject to}& & \sum_{j=1}^{n_n}z_j =m\\
    &&& Gz \leq b\\
    &&& u_{ij} \leq z_j,\quad i \in \{1,\ldots,n_n\},\; j \in \{1,\ldots,n_n\}\\
    &&&\sum_{j=1}^{n_n}u_{ij} = 1, \quad  i \in \{1,\ldots,n_n\}\\
    &&& u_{ii} = 0, \quad  i \in \{1,\ldots,n_n\}\\
    &&& u_{ij} \geq 0 \quad i \in \{1,\ldots,n_n\},\; j \in \{1,\ldots,n_n\}\\
    &&& z \in [0,1]^{n_n}.
    \end{alignedat}
\end{equation}
In comparison, when $f(z) = f_D(z)$, Problem~\eqref{eq:prob_convex} yields:
\begin{equation}
    \label{eq:dopt_convex_0}
    \begin{alignedat}{3}
    &\text{minimize}&\; \; &-\log\,\det\bigg(X_0+\sum_{j=1}^{n_n}z_jA_jA_j^T\bigg)\\
    &\text{subject to}& & \sum_{j=1}^{n_n}z_j=m\\
    &&&Gz \leq b\\
    &&& z \in [0,1]^{n_n},
    \end{alignedat}
\end{equation}
which is equivalently re-written as a semidefinite program:
\begin{equation}
    \label{eq:dopt_convex}
    \begin{alignedat}{3}
    &\text{minimize}&\; \; &-\sum_{k=1}^{n_r}y_k\\
    &\text{subject to}& & (y_k,1,W_{kk}) \in \mathcal{E},  \quad k \in \{1,\ldots,n_r\}\\
    &&& \begin{bmatrix}X_0+\sum_{j=1}^{n_n} z_jA_jA_j^T & W\\
        W^T & \diag(W)\end{bmatrix} \succeq 0\\
    &&& W \text{ is lower triangular}\\
    &&&\sum_{j=1}^{n_n}z_j=m\\
    &&& Gz \leq b\\
    &&& z \in [0,1]^{n_n},
    \end{alignedat}
\end{equation}
where $\text{diag}(W) \in \mathbb{R}^{n_r \times n_r}$ is the diagonal matrix with the same diagonal elements as $W$, and $\mathcal{E}$ is the exponential cone. Problem~\eqref{eq:dopt_convex} is a semidefinite program with $n_r+n_n$ scalar variables, one matrix variable, and $n_r$ exponential cone constraints. Finally, when $f(z) = f_{\beta}(z)$, Problem~\eqref{eq:prob_convex} is:
\begin{equation}
    \label{eq:cheb_convex}
    \begin{alignedat}{3}
        &\text{minimize}&\; \; & t\\
        &\text{subject to}& & f_D(z) \leq \frac{t}{w^{\beta}_D} + f_D^*\\
        &&& f_T(z) \leq \frac{t}{w^{\beta}_T} + f_T^*\\
        &&& \sum_{j=1}^{n_n}z_j=m\\
        &&& Gz \leq b\\
        &&& z \in [0,1]^{n_n}.
    \end{alignedat}
\end{equation}
By Lemma~\ref{lem:pmedian_equivalently}, Problem~\eqref{eq:cheb_convex} is equivalent to:
\begin{equation}
    \label{eq:cheb_convex2}
    \begin{alignedat}{3}
        &\text{minimize}&\; \; & t\\
        &\text{subject to}& & -\log\,\det\bigg(X_0+\sum_{j=1}^{n_n}z_jA_jA_j^T\bigg) \leq \frac{t}{w^{\beta}_D} + f_D^*\\
        &&& \sum_{i=1}^{n_n}\sum_{j=1}^{n_n}C_{ij}u_{ij} \leq \frac{t}{w^{\beta}_T} + f_T^*\\
        &&& u_{ij} \leq z_j,\quad i \in \{1,\ldots,n_n\},\; j \in \{1,\ldots,n_n\}\\
        &&&\sum_{j=1}^{n_n}u_{ij} = 1, \quad  i \in \{1,\ldots,n_n\}\\
        &&& u_{ii} = 0, \quad  i \in \{1,\ldots,n_n\}\\
        &&& u_{ij} \geq 0 \quad i \in \{1,\ldots,n_n\},\; j \in \{1,\ldots,n_n\}\\
        &&& \sum_{j=1}^{n_n}z_j=m\\
        &&&Gz \leq b\\
        &&& z \in [0,1]^{n_n},
    \end{alignedat}
\end{equation}
which is re-written into the following semidefinite program:
\begin{equation}
    \label{eq:cheb_convex3}
    \begin{alignedat}{3}
    &\text{minimize}&\; \; & t\\
    &\text{subject to}& &  -\sum_{k=1}^{n_r}y_k \leq \frac{t}{w^{\beta}_D} + f^*_D\\
    &&& \sum_{i=1}^{n_n}\sum_{j=1}^{n_n} C_{ij}u_{ij}\leq \frac{t}{w^{\beta}_T} + f^*_T \\
    &&& (y_k,1,W_{kk}) \in \mathcal{E},  \quad k \in \{1,\ldots,n_r\}\\
    &&& \begin{bmatrix}X_0+\sum_{j=1}^{n_n} z_jA_jA_j^T & W\\
        W^T & \diag(W)\end{bmatrix} \succeq 0\\
    &&& W \text{ is lower triangular}\\
    &&& u_{ij} \leq z_j,\quad i \in \{1,\ldots,n_n\},\; j \in \{1,\ldots,n_n\}\\
    &&&\sum_{j=1}^{n_n}u_{ij} = 1, \quad  i \in \{1,\ldots,n_n\}\\
    &&& u_{ii} = 0, \quad  i \in \{1,\ldots,n_n\}\\
    &&& u_{ij} \geq 0 \quad i \in \{1,\ldots,n_n\},\; j \in \{1,\ldots,n_n\}\\
    &&& \sum_{j=1}^{n_n}z_j =m\\
    &&& Gz \leq b\\
    &&& z \in [0,1]^{n_n}.
    \end{alignedat}
\end{equation}

\section*{Acknowledgements}
\noindent This work was supported by EPSRC (EP/P004229/1, Dynamically Adaptive and Resilient Water Supply Networks for a Sustainable Future).

\bibliography{references}

\end{document}